\documentclass{amsart}

\usepackage{amssymb}		
\usepackage{eucal}		

\def\varinfty{\underline{\infty}}

\def\propto{\varpropto}
\renewcommand\ge\geqslant
\renewcommand\le\leqslant

\newtheorem{theorem}{Theorem}[section]
\newtheorem{corollary}[theorem]{Corollary}
\newtheorem{proposition}[theorem]{Proposition}
\newtheorem{lemma}[theorem]{Lemma}

\theoremstyle{definition} 
\newtheorem{definition}[theorem]{Definition}

\theoremstyle{remark}

\newtheorem{remarkaftertheorem}{Remark}[theorem]

\numberwithin{equation}{section}

\newcommand{\Hl}{\mathop{{}\it Hl}\nolimits}

\begin{document}

\title[Algebraic Hypergeometric Transformations of Modular Origin]{Algebraic Hypergeometric Transformations\\ of Modular Origin}

\author{Robert S. Maier}
\address{Depts.\ of Mathematics and Physics, University of Arizona, Tucson,
AZ 85721, USA}
\email{rsm@math.arizona.edu}
\urladdr{http://www.math.arizona.edu/\~{}rsm}
\thanks{The author was supported in part by NSF Grant No.\ PHY-0099484.}

\renewcommand{\subjclassname}{\textup{2000} Mathematics Subject Classification}

\subjclass{Primary 11F03; 11F20, 33C05.}
\date{}

\begin{abstract}
It is shown that Ramanujan's cubic transformation of the Gauss
hypergeometric function~${}_2F_1$ arises from a relation between modular
curves, namely the covering of~$X_0(3)$ by~$X_0(9)$.  In~general,
when~$2\le N\le 7$ the $N$-fold cover of~$X_0(N)$ by~$X_0(N^2)$ gives rise
to an algebraic hypergeometric transformation.  The $N=2,3,4$
transformations are arithmetic--geometric mean iterations, but the
$N=5,6,7$ transformations are new.  In~the final two the change of
variables is not parametrized by rational functions, since $X_0(6),X_0(7)$
are of genus~$1$.  Since their quotients $X_0^+(6),X_0^+(7)$ under the
Fricke involution (an~Atkin--Lehner involution) are of genus~$0$, the
parametrization is by two-valued algebraic functions.  The resulting
hypergeometric transformations are closely related to the two-valued
modular equations of Fricke and H.~Cohn.
\end{abstract}

\maketitle

\section{Introduction}

Identifiable functions that satisfy functional equations of high degree are
rare flowers.  For this reason, much attention has been paid to Ramanujan's
parametrized cubic transformation
\begin{equation}
\label{eq:Rama}
{}_2F_1\left(\tfrac13,\tfrac23;\,1;\,1-\Bigl(\frac{1-x}{1+2x}\Bigr)^3\right)
=
(1+2x)\,\, {}_2F_1\left(\tfrac13,\tfrac23;\,1;\,x^3\right)
\end{equation}
of a particular case of the Gau\ss\ hypergeometric
function~${}_2F_1(\alpha,\beta;\gamma;\cdot)$, namely
${}_2F_1(\frac13,\frac23;1;\cdot)$.  In Ramanujan's theory of elliptic
integrals in signature~$3$, the functional equation~(\ref{eq:Rama}) appears
as the cubic arithmetic--geometric mean iteration~\cite{Berndt95}.  The
first published proof, by the Borweins~\cite{Borwein91}, relied on the
existence of an action of the modular group~$\Gamma$ on the right-hand
argument~$x^3$, and the invariance of that argument under a subgroup.
Additional proofs have appeared~\cite{Borwein94,Chan98}.

In this article we show that (\ref{eq:Rama})~is one of six related
hypergeometric identities, naturally indexed by an integer~$N=2,\dots,7$.
The $N$'th identity arises from the modular curve $X_0(N)$, in~particular
from its $N$-sheeted covering by~$X_0(N^2)$.  Here the curve~$X_0(N)$,
which classifies $N$-isogenies of elliptic curves over~$\mathbb{C}$, is the
(compactified) quotient of the upper half-plane~$\mathfrak{H}$ by the Hecke
congruence subgroup $\Gamma_0(N)$ of~$\Gamma$.  The connection to the
Borweins' proof should be evident.  Of~the six identities, the
identity~(\ref{eq:Rama}) is associated with~$X_0(3)$, and the ones
associated with $X_0(2)$ and~$X_0(4)$ were also found by Ramanujan.  (They
are the quadratic iteration in signature~$4$ and the quartic one in
signature~$2$, respectively.)  The identities associated with $X_0(N)$,
$N=5,6,7$, are new.  The reason for their not having been discovered
previously may be that they are most briefly expressed not in~terms
of~${}_2F_1$, but rather in~terms of~$\Hl$, the so-called local Heun
function~\cite{Ronveaux95}.  The functions ${}_2F_1,\Hl$ are solutions of
canonical Fuchsian differential equations on~$\mathbb{P}^1(\mathbb{C})$
with three and four singular points, respectively.

Each identity is a modular equation.  When $N=2,3,4,5$, the curve
$X_0(N^2)$~is of genus~$0$ and has a Hauptmodul (global uniformizing
parameter), in~terms of which the identity can be expressed.  When $N=6,7$,
the curve $X_0(N^2)$ is of genus~$1$ but its quotient $X_0^+(N^2)$ under the
Fricke involution has a Hauptmodul.  So the $N=6,7$ identities can be
formulated without elliptic functions by including algebraic constraints.
All of these hypergeometric identities can be viewed as algebraic
transformations of series, but the following statement of the quintic
($N=5$) identity in series language makes it clear how difficult it would
be to prove them by series manipulation.  The accompanying proof is really
a {\em verification\/}, and will be replaced in the sequel by a derivation
based on Hauptmoduln.  The function~$h_5$ can be expressed in~terms of
$\Hl$ or~${}_2F_1$.  It~will be shown that
$h_5(z)=[\frac15(z^2+10z+5)]^{-1/4}\,{}_2F_1\left(\frac1{12},\frac5{12};1;1728z/(z^2+10z+5)^3\right)$.

\begin{proposition}
\label{prop:h5}
Let $h_5$, a $\mathbb{C}$-valued function, be defined in a neighborhood of\/
$0\in\mathbb{C}$ by $h_5(z)=\sum_{n=0}^\infty c_nz^n$, where the
coefficients satisfy the three-term recurrence
\begin{equation}
\label{eq:h5recurrence}
(2n-1)^2\,c_{n-1} + 2(44n^2+22n+5)\,c_n + 500(n+1)^2\,c_{n+1}=0,
\end{equation}
initialized by $c_{-1}=0$ and $c_0=1$.  Then for all~$x$ in a neighborhood
of\/~$0$, 
\begin{multline}
\label{eq:quintic}
h_5\left(x(x^4+5x^3+15x^2+25x+25)\right)\\
= 5\left[x^4+5x^3+15x^2+25x+25\right]^{-1/2} h_5\left(\frac{x^5}{x^4+5x^3+15x^2+25x+25}\right).
\end{multline}
That is, $h_5$~satisfies a quintic functional equation.
\end{proposition}
\begin{remarkaftertheorem}
The sequence $d_n:=500^nc_n$, $n\ge0$, of Maclaurin coefficients of
$h_5(500z)$ is an integral sequence.  It~begins $1$, $-10$, $230$, $-6500$,
$199750$, $-6366060$, $204990300$, $-6539387400$, $\dots$.
\end{remarkaftertheorem}
\begin{proof}[Proof of Proposition\/~{\rm\ref{prop:h5}}]
In a neighborhood of $z=0$, $h_5=h_5(z)$ is analytic and satisfies the
second-order differential equation
\begin{equation}
\label{eq:h5de}
\left\{D_z^2 + \left[\frac1z+\frac{z+11}{z^2+22z+125}\right]D_z + 
\left[
\frac{z+10}{4z(z^2+22z+125)}
\right]\right\}h_5=0,
\end{equation}
as can be verified by termwise differentiation of its defining series.  By
changing variables in~(\ref{eq:h5de}), it~can be shown that as functions
of~$x$, the two sides of~(\ref{eq:quintic}) satisfy a common differential
equation, namely
\begin{equation}
\label{eq:h5de5}
\left\{D_x^2 + \left[\frac1x + \frac{D_xA(x)}{A(x)} + \frac{D_xB(x)}{2B(x)} \right]D_x
+\left[\frac{25(x\,A(x)+10)}{4x\,A(x)\,B(x)}\right]\right\}f = 0,
\end{equation}
with $A(x):=x^4+5x^3+15x^2+25x+25$ and $B(x):=x^2+2x+5$.  The point $x=0$
is a regular singular point of~(\ref{eq:h5de5}), with characteristic
exponents~$0,0$.  It~follows from the theory of Fuchsian differential
equations that in a neighborhood of~$x=0$, there is a unique analytic
solution of~(\ref{eq:h5de5}) that equals unity at~$x=0$.  But both sides
of~(\ref{eq:quintic}) are analytic at~$x=0$ and equal unity there; so each
can be identified with this unique solution, and they must equal each
other.
\end{proof}

This article is organized as follows.  Basic facts about modular curves and
their coverings are summarized in Section~\ref{sec:background}.  In
Sections~\ref{sec:liftingsofcusps} and~\ref{sec:liftingsofDEs} the covering
of $X_0(N)$ by~$X_0(N^2)$ is reviewed, and useful lemmas are proved.
Section~\ref{sec:keyresults} contains the key results, Theorems~\ref{thm:1}
and~\ref{thm:2}, which apply to $N=2,3,4,5$ and~$N=6,7$ respectively.  The
algebraic hypergeometric transformations derived from these theorems are
worked~out in Sections~\ref{subsec:g0} and~\ref{subsec:g1}.  The derivation
relies on a number of formulas relating canonical Hauptmoduln, which are
collected in an appendix.  In Section~\ref{sec:extensions} some possible
extensions are indicated.

\section{Curves and Coverings}
\label{sec:background}

The left action of ${\it SL}(2,\mathbb{Z})$
on~$\left(\begin{smallmatrix}\tau_1\\\tau_2\end{smallmatrix}\right)$ with
$\tau_1,\tau_2\in\mathbb{C}^*$ independent over~$\mathbb{R}$, i.e., on
elliptic curves of the form $\mathbb{C}/\langle\tau_1,\tau_2\rangle$,
induces a projective action of the modular group $\Gamma:={\it
SL}(2,\mathbb{Z})/\{\pm I\}$ on the upper half plane
$\mathfrak{H}\ni\tau:=\tau_1/\tau_2$.  The quotient
$\Gamma\setminus\mathfrak{H}$ is the space of isomorphism classes of
elliptic curves over~$\mathbb{C}$, and its compactification
$\Gamma\setminus\left(\mathfrak{H}^*:=\mathfrak{H}\cup\mathbb{Q}\cup\{{\rm
i}\infty\}\right)$, with cusps, is the modular curve~$X(1)$.

The isomorphism classes of $N$-isogenies of elliptic curves $\phi:E\to E'$
with ${N>1}$, i.e., isogenies with kernel equal to the cyclic group~$C_N$,
are similarly classified by the non-cusp points of the modular curve
$X_0(N):=\Gamma_0(N)\setminus\mathfrak{H}^*$, where
$\Gamma_0(N)=\left\{\left(\begin{smallmatrix}a&b\\c&d\end{smallmatrix}\right)\in\Gamma\mid
N|c\right\}$ is the $N$'th congruence subgroup.  For any such $N$-isogeny
of elliptic curves $\phi:E\to E'$, i.e.,
$\phi:\mathbb{C}/\langle\tau,1\rangle\to\mathbb{C}/\langle\tau',1\rangle$,
there is a dual one
$\phi^*:\mathbb{C}/\langle\tau',1\rangle\to\mathbb{C}/\langle\tau,1\rangle$.
This correspondence yields the Fricke involution $w_N$ on~$X_0(N)$, defined
as~$E\leftrightarrow E'$.  On~the level of the unquotiented half-plane
$\mathfrak{H}\ni\tau$ it~is simply the map $\tau\mapsto-1/N\tau$.  The
curve~$X_0^+(N)$ is defined as the quotient of~$X_0(N)$ by~$\langle
w_N\rangle$, the group of two elements generated by~$w_N$.  Each non-cusp
point of~$X_0^+(N)$ corresponds to an unordered pair $\bigl\{\phi:E\to E',\
\phi^*:E'\to E\bigr\}$.

The curve $X(1)$ is of genus zero, so its function field is singly
generated; the generator can be taken to be the Klein--Weber $j$-invariant.
This is effectively a meromorphic function on~$\mathfrak{H}^*$ with an
expansion about~$\tau={\rm i}\infty$ that begins $q^{-1}+\cdots$, where
$q:=e^{2\pi{\rm i}\tau}$~is the local uniformizing parameter.  For
any~$N>1$, the function field of~$X_0(N)$ is $\mathbb{C}(j,j_N)$, where
$j_N(\tau):=j(N\tau)$.  The Fricke involution~$w_N$ interchanges~$j,j_N$.
If $g(X_0(N))=0$ then $X_0(N)$ will have a Hauptmodul, which may be
denoted~$x_N=x_N(\tau)$, and both $j$~and~$j_N$ will necessarily be
rational functions of~it.  Being univalent, $x_N$~can be chosen to have a
simple zero (resp.~pole) at the cusp~$\tau={\rm i}\infty$ (resp.~$\tau=0$).
Since these two cusps are interchanged by~$w_N$, $x_N$~can be defined so
that $x_N(\tau)\,x_N(-1/N\tau)=\kappa_N$ for any specified
$\kappa_N\in\mathbb{Q}^*$.

If~$N|N'$ then $\Gamma_0(N)\ge \Gamma_0(N')$, yielding a (ramified)
covering of~$X_0(N)$ by $X_0(N')$ and an injection of the function field
of~$X_0(N)$ into that of~$X_0(N')$.  If, for example, both have
Hauptmoduln, then $x_N$~will be a rational function of~$x_{N'}$.  If
$g(X_0(N'))>0$ but $g(X_0^+(N'))=0$, then $X_0(N')$~will be a hyperelliptic
curve that doubly covers~$X_0^+(N')$, and $X_0^+(N')$~will have a
Hauptmodul~$t_{N'}$ such that $X_0(N')$~is defined by $s_{N'}^2={\sf
P}(t_{N'})$ for some polynomial~$\sf P$ of degree $2g(X_0(N'))+2$.  The
hyperelliptic involution $s_{N'}\mapsto -s_{N'}$ will be the Fricke
involution~$w_{N'}$.  The function field of~$X_0(N')$ will be generated
by~$t_{N'},s_{N'}$; and like~$s_{N'}$, the Hauptmodul~$x_N$ of~$X_0(N)$
will be a two-valued algebraic function of~$t_{N'}$.

This article will focus on the especially interesting case $N'=N^2$.  When
$N=2,3,4$, or~$6$, the automorphism group~$\Gamma_0(N^2)$ is conjugate
in~${\it PSL}(2,\mathbb{R})$ to the level-$N$ principal congruence
subgroup~$\Gamma(N)$ of~$\Gamma$, so the corresponding curves
$X_0(N^2),X(N)$ are isomorphic; and for any~$p$, the quotient
curve~$X_0^+(p^2)$ is isomorphic to the arithmetically important
curve~$X_{\rm split}(p)$.  Computational treatment of the covering
of~$X_0(N)$ by~$X_0(N^2)$ is facilitated by the many known formulas
relating the associated Hauptmoduln, originating largely with Fricke.
These do~not appear in~full in any modern reference, so they are reproduced
here in an appendix, in an enhanced format.  Each~$N>1$ with $g(X_0(N))=0$
and either $g(X_0(N^2))=0$ or $g(X_0^+(N^2))=0$ is included.  The values
of~$N$ turn~out to be $N=2,\dots,7$.

The genus $g(X_0(N))$ here comes from the Hurwitz formula, or directly from
Euler's theorem~\cite{Schoeneberg74}.  The covering $j:X_0(N)\to
\mathbb{P}^1(\mathbb{C})\cong X(1)$ is $\psi(N)$-sheeted, where
$\psi(N):=N\prod_{p|N}(1+1/p)$ is the index $[\Gamma:\Gamma_0(N)]$.  It is
ramified only above the cusp~$j=\infty$ and the elliptic fixed points
${j=0,12^3}$, corresponding to equianharmonic and lemniscatic elliptic
curves respectively; i.e., only above the vertices $\tau={\rm
i}\infty,\rho:=e^{2\pi i/3},{\rm i}$ of the fundamental half-domain
of~$\Gamma$ in~$\mathfrak{H}^*$.  The fibre above~$j=\infty$ consists of
$\sigma_\infty(N) := \sum_{d|N}\varphi\left((d,N/d)\right)$ cusps, where
$\varphi$~is the Euler totient function, i.e.,
$\varphi(N):=N\prod_{p|N}(1-1/p)$.  The fibre above $j=0$ (resp.\ $j=12^3$)
includes $\varepsilon_\rho(N)$ cubic elliptic points (resp.\
$\varepsilon_{\rm i}(N)$ quadratic ones), each with unit multiplicity;
other points, if~any, have cubic (resp.\ quadratic) multiplicity.  Here
\begin{displaymath}
\varepsilon_\rho(N):=
\begin{cases}
\prod_{p|N}\left(1+\left(\frac{-3}{p}\right)\right), & 9\!\!\not|N, \\
0, & 9|N,
\end{cases}
\qquad
\varepsilon_{\rm i}(N):=
\begin{cases}
\prod_{p|N}\left(1+\left(\frac{-1}{p}\right)\right), & 4\!\!\not|N, \\
0, & 4|N,
\end{cases}
\end{displaymath}
with $\left(\frac{\cdot}{\cdot}\right)$ the Legendre symbol.  An
application of Euler's theorem yields
\begin{equation}
g\left(X_0(N)\right) = 1+\frac{\psi(N)}{12} - \frac{\sigma_\infty(N)}2 -
\frac{\varepsilon_\rho(N)}{3} - \frac{\varepsilon_{\rm i}(N)}{4}.
\end{equation}
A further computation, first performed by Fricke~\cite{Fricke22,Kluit77},
reveals that in
\begin{equation}
g\left(X_0^+(N)\right)=\tfrac12\left[g\left(X_0(N)\right)+1\right]-\tfrac14a(N),
\end{equation}
which follows from the Hurwitz formula, with $a(N)$~the number of fixed
points of the Fricke involution on~$X_0(N)$, the quantity~$a(N)$
(when~$N\ge5$, at~least) equals $h(-4N)+h(-N)$ if~$N\equiv1\pmod4$ and
$h(-4N)$ otherwise.  Here $h(-d)$~is the class number of the imaginary
quadratic field~$\mathbb{Q}(\sqrt{-d})$.

\section{Liftings of Cusps}
\label{sec:liftingsofcusps}

The cusps of~$X_0(N)$ have the following
description~\cite{Gonzalez91,Ogg73}.  The set of cusps of~$\mathfrak{H}^*$,
$\mathbb{P}^1(\mathbb{Q})=\mathbb{Q}\cup\{{\rm i}\infty\}\ni\tau$, is
partitioned into classes, each equivalent under~$\Gamma_0(N)$.  A~system of
representatives, i.e., a choice of one from each class, may be taken to
comprise certain fractions $\tau=\frac{a}{d}$ for each~$d|N$, with $1\le
a<d$ and~$(a,d)=1$.  Here $a$~is reduced modulo $f_{d,N}:=(d,N/d)$ while
remaining coprime to~$d$, so there are $\varphi((d,N/d))$ values of~$a$,
and hence $\varphi((d,N/d))$ cusps in~$\mathbb{P}^1(\mathbb{Q})$ associated
to~$d$, which are inequivalent under~$\Gamma_0(N)$.  This is the source of
the above formula for~$\sigma_\infty(N)$.  Each cusp of the
form~$\frac{a}d$ has width $e_{d,N}:=N/df_{d,N}$, i.e., multiplicity
$e_{d,N}$ above~$X(1)$.  That~is, the fibre of the $\psi(N)$-sheeted cover
$\pi_N:X_0(N)\to X(1)$ above the unique cusp of~$X(1)$, located
at~$j=\infty$, includes (the equivalence class~of) $\tau=\frac{a}d$ with
multiplicity~$e_{d,N}$.  To~emphasize that a cusp of~$X_0(N)$ is an
equivalence class, the notation~$\left[\frac{a}d\right]$ will be used,
or~$\left[\frac{a}d\right]_N$ if the modular curve needs to be indicated.
Since the `distinguished' cusps $\tau=0,{\rm i}\infty$ of~$X_0(N)$ are
equivalent to~$\tau=\frac{1}1,\frac{1}N$ respectively, they may be written
as~$\left[\frac11\right]_N,\left[\frac1N\right]_N$.

Now consider the inverse images of the $\sigma_\infty(N)$ cusps of~$X_0(N)$
under its $N$-sheeted covering by~$X_0(N^2)$, which will be
denoted~$\phi_N$.  Note first that the covering $\pi_{N^2}:X_0(N^2)\to
X(1)$, which has $\psi(N^2)=N\psi(N)$ sheets, satisfies
$\pi_{N^2}=\pi_N\circ\phi_N$.  The coverings
$X_0(N^2)\stackrel{\phi_N}{\longrightarrow}X_0(N)\stackrel{\pi_N}{\longrightarrow}X(1)$
correspond to the subgroup relations $\Gamma_0(N^2)<\Gamma_0(N)<\Gamma$.
For any~$N$, $\Gamma_0(N)$~is normalized in~${\it PSL}(2,\mathbb{R})$ by
the Fricke involution $w_N:\tau\mapsto-1/N\tau$.  That~is, if
$\left(\begin{smallmatrix}a&b\\c&d\end{smallmatrix}\right)\in\Gamma_0(N)$
then
$\left(\begin{smallmatrix}0&-1\\N&0\end{smallmatrix}\right)^{-1}\left(\begin{smallmatrix}a&b\\c&d\end{smallmatrix}\right)\left(\begin{smallmatrix}0&-1\\N&0\end{smallmatrix}\right)\in\Gamma_0(N)$;
and $w_N$~induces a permutation of the cusps of~$X_0(N)$.  This
interchanges any cusp $\left[\frac{a}d\right]$ with some cusp
$\bigl[\frac{a'}{d'}\bigr]$, where~$d'=N/d$.

There is another covering by~$X_0(N^2)$ that will be equally important.
Though $\Gamma_0(N^2)$ is normalized by~$w_{N^2}$, $\Gamma_0(N)$~is~not:
it~is conjugated to an isomorphic subgroup
$\Gamma_0(N)':=w_{N^2}^{-1}\Gamma_0(N)w_{N^2}$ of~${\it
PSL}(2,\mathbb{R})$.  Since $\Gamma_0(N^2)<\Gamma_0(N)'$, there is a
corresponding cover $X_0(N^2)\stackrel{\phi_N'}{\longrightarrow}X_0(N)'$,
where $X_0(N)':=\Gamma_0(N)'\setminus\mathfrak{H}^*$.  Since
$\Gamma_0(N)'$, unlike~$\Gamma_0(N)$, is not a subgroup of the modular
group~$\Gamma$, the curve $X_0(N)'$ does~not naturally cover~$X(1)$.  But
if $X_0(N)$ has a Hauptmodul~$x_N$, which may be regarded as a
$\Gamma_0(N)$-invariant function $x_N(\tau)$ on~$\mathfrak{H}^*$, then
$X_0(N)'$ will too, and it may be chosen to be the $X_0(N)'$-invariant
function $x'_N(\tau):=x_N(N\tau)$.  This is because
$\Gamma_0(N)'=w^{-1}\Gamma_0(N)w$ where
$w=\left(\begin{smallmatrix}N&0\\0&1\end{smallmatrix}\right)$, i.e., $w$~is
the map~$\tau\mapsto N\tau$.

The cusps of $X_0(N),X_0(N)'$, which are subsets
of~$\mathbb{P}^1(\mathbb{Q})$ invariant under $\Gamma_0(N)$,
$\Gamma_0(N)'$, are related thus: $[\tau]_N$~is a cusp of~$X_0(N)$ iff
$[\tau]_N/N:=\{\,x/N\mid x\in[\tau]_N\,\}$~is a cusp of~$X_0(N)'$.  The
following lemma specifies how the two most important cusps of~$X_0(N)$
(resp.~$X_0(N)'$) are lifted to formal sums of cusps of~$X_0(N^2)$.
\begin{lemma}
\label{lem:cusplifting}
Inverse images of the distinguished cusps $\tau=0,{\rm i}\infty$
on~$X_0(N)$ {\rm(}i.e., of\/
$\left[\frac11\right]_N,\left[\frac1N\right]_N${\rm)} under the $N$-sheeted
cover $\phi_N:X_0(N^2)\to X_0(N)$ are given by
\begin{displaymath}
{\phi_N}^{-1}(\left[\tfrac11\right]_N) = N\cdot\left[\tfrac11\right]_{N^2}, \qquad\qquad
{\phi_N}^{-1}(\left[\tfrac1N\right]_N) = \!\sum_{\substack{d\ {\rm s.t.}\,N|d|N^2,\\a=a_1,\dots,a_{\varphi(f_{d,N^2})}}}
\!\!\left[\tfrac{a}d\right]_{N^2},
\end{displaymath}
where the right-hand sides list cusps of~$X_0(N^2)$, and premultiplication
by~$N$ indicates multiplicity.  The sum is over\/ {\rm(}i\,{\rm)} $d$~such
that $d|N^2$ with~$N|d$, and\/ {\rm(}ii\,{\rm)}~the corresponding
$\varphi(f_{d,N^2})$ values of~$a$, giving $N$~terms in~all.  Inverse
images of the cusps $\tau=0,{\rm i}\infty$ on~$X_0(N)'$ under
$\phi_N':X_0(N^2)\to X_0(N)'$ are similarly given by
\begin{displaymath}
{\phi_N'}^{-1}(\left[\tfrac11\right]_N/N) = \!\sum_{\substack{d\ {\rm s.t.}\,d|N|N^2,\\a=a_1,\dots,a_{\varphi(f_{d,N^2})}}}
\!\!\left[\tfrac{a}d\right]_{N^2}, \qquad\qquad
{\phi_N'}^{-1}(\left[\tfrac1N\right]_N/N) = N\cdot\left[\tfrac1{N^2}\right]_{N^2}.
\end{displaymath}
The sum is over $d$~such that~$d|N$, and over the corresponding
$\varphi(f_{d,N^2})$ values of~$a$.
\end{lemma}
\begin{proof}
For the duration of this proof, write any cusp~$\left[\tfrac{a}d\right]$
as~$\left(\begin{smallmatrix}a\\d\end{smallmatrix}\right)\in\mathbb{Z}^2$,
to permit left-multiplication by elements of~${\it SL}(2,\mathbb{Z})$.  The
cosets of~$\Gamma_0(N)$ in~$\Gamma_0(N^2)$ are represented
by~$\left(\begin{smallmatrix}1&0\\kN&1\end{smallmatrix}\right)$,
$k=0,\dots,N-1$.  The
cusp~$\left(\begin{smallmatrix}1\\1\end{smallmatrix}\right)$ of~$X_0(N)$ is
lifted by~$\phi_N$ to the $N$~cusps
$\mathfrak{c}_k:=\left(\begin{smallmatrix}1&0\\kN&1\end{smallmatrix}\right)
\left(\begin{smallmatrix}1\\1\end{smallmatrix}\right)
=\left(\begin{smallmatrix}1\\kN+1\end{smallmatrix}\right)$,
$k=0,\dots,N-1$.  To~prove they are the same in~$X_0(N^2)$, it~suffices to
find $\gamma_k\in\Gamma_0(N^2)$ such that
$\gamma_k\cdot\mathfrak{c}_k=\mathfrak{c}_0$.  By examination,
$\left(\begin{smallmatrix}k^2N^2+kN+1&-kN\\k^2N^2&-kN+1\end{smallmatrix}\right)$
will work.

The lifting of $\left(\begin{smallmatrix}1\\N\end{smallmatrix}\right)$
similarly comprises the $N$~cusps
$\mathfrak{d}_k:=\left(\begin{smallmatrix}1&0\\kN&1\end{smallmatrix}\right)
\left(\begin{smallmatrix}1\\N\end{smallmatrix}\right)
=\left(\begin{smallmatrix}1\\(k+1)N\end{smallmatrix}\right)$,
$k=0,\dots,N-1$, of~$X_0(N^2)$.  Since $w_{N^2}\mathfrak{d}_k=
\left(\begin{smallmatrix}-(k+1)N\\N^2\end{smallmatrix}\right) \sim
\left(\begin{smallmatrix}-(k+1)\\N\end{smallmatrix}\right)$, the cusps
$\{w_{N^2}\mathfrak{d}_k\}_{k=0}^{N-1}$ are of the
form~$\left(\begin{smallmatrix}a\\d\end{smallmatrix}\right)$, where
$d$~runs over the divisors of~$N$ and $a$~runs over the
$\varphi(d)=\varphi((d,N^2/d))$ integers in the range $1,\dots,d-1$ that
are relatively prime to~$d$.  As noted, $w_{N^2}$~maps any cusp
$\left(\begin{smallmatrix}a\\d\end{smallmatrix}\right)$ of~$\Gamma_0(N^2)$
to a cusp $\left(\begin{smallmatrix}a'\\d'\end{smallmatrix}\right)$ with
$d'=N^2/d$; so applying~$w_{N^2}$ again yields the given formula for the
formal sum~$\sum_{k=0}^{N-1}\mathfrak{d}_k$.

The statements about liftings by~$\phi'_N$ can be proved likewise; or
simply by applying~$w_{N^2}$ to both sides of each of the two previously
derived formulas.
\end{proof}

The following lemma subsumes the half of Lemma~\ref{lem:cusplifting} that
deals with~$\phi_N$.  It~can be proved by a similar argument.
\begin{lemma}
\label{lem:cusplifting2}
The fibre of\/ $\phi_N$ over any cusp of~$X_0(N)$ of the
form\/~$\left[\frac{a}d\right]_N$ with~$d|N$ consists of $d/(d,N/d)$ cusps
of~$X_0(N^2)$, each of the form\/ $\left[\frac{\tilde a}{\tilde
d}\right]_{N^2}$ with $\tilde d|N^2$, where $\tilde d$~is constrained to
satisfy $d=(\tilde d,N)$.  Each of these cusps of~$X_0(N^2)$ appears in the
fibre with multiplicity equal to the width quotient
$e_{d,N^2}/e_{d,N}=N(d,N/d)/d$.
\end{lemma}

A consequence of the lemma is that for any~$d|N$, the inverse image
under~$\phi_N$ of the set consisting of the $\varphi((d,N/d))$ cusps
of~$X_0(N)$ of the form~$\left[\frac{a}d\right]_N$ is a set consisting of
the $[d/(d,N/d)]\,\varphi((d,N/d))$ cusps of~$X_0(N^2)$ of the
form~$\left[\frac{\tilde a}{\tilde d}\right]_{N^2}$, where $\tilde
d$~ranges over the solutions of~$d=(\tilde d,N)$ with~$\tilde d|N^2$.
Each of these cusps of~$X_0(N^2)$ appears in the inverse image with
multiplicity~$N(d,N/d)/d$.

\section{Liftings of Differential Operators}
\label{sec:liftingsofDEs}

\subsection{Some background}

The coordinate~$\tau$ on~$\mathfrak{H}^*$ can be viewed as a multivalued
function on any algebraic curve~$X_1$ of the form
$\Gamma_1\setminus\mathfrak{H}^*$, where $\Gamma_1$~is a Fuchsian subgroup
(of~the first kind) of the automorphism group~${\it PSL}(2,\mathbb{R})$.
It~is well~known that on~$X_1$, $\tau$~will locally equal the ratio of two
independent solutions of an appropriate second-order Fuchsian differential
equation~\cite{Ford51,Harnad2000}.  For example, consider the case
$\Gamma_1=\Gamma$, since $X_1=X(1)$ is parametrized by the $j$-invariant
and the equation may be written in~terms of the derivation $D_j:=d/dj$.
Actually, it~is more convenient to use the Hauptmodul $\hat J:=12^3/j$, the
reciprocal of Klein's invariant~$J=j/12^3$.
\begin{proposition}
\label{prop:0}
In a neighborhood of any point
on~$X(1)=\Gamma\setminus\mathfrak{H}^*\ni\hat J$, any branch of~$\tau$
equals the ratio $u_1/u_2$ of two independent local solutions of
\begin{equation}
\label{eq:lowest}
L_{\frac1{12},\frac5{12};1}\, u:=
\left\{D_{\!\hat J}^2
+ \left[ \frac1{\hat J} + \frac1{2(\hat J-1)} \right] D_{\!\hat J}
+ \frac{5/144}{\hat J(\hat J-1)}\right\} u = 0.
\end{equation}
Moreover, any such ratio is of the form $(a\tau+b)/(c\tau+d)$ with
$ad\neq bc$.
\end{proposition}
Equation~(\ref{eq:lowest}) is a Gau\ss\ hypergeometric equation on the $\hat
J$-sphere~$X(1)$.  It~is Fuchsian, with regular singular points at~$\hat
J=0,1,\infty$ and respective characteristic exponents
$0,0;\,0,\frac12;\,\frac{1}{12},\frac{5}{12}$.  It will be referred to as
${}_2\mathcal{E}_1(\frac1{12},\frac5{12};1)$, since the unique solution
analytic at~$\hat J=0$ and normalized to unity there is the hypergeometric
function ${}_2F_1(\frac1{12},\frac5{12};1;\cdot)$.
On~${\mathfrak{H}}^*\ni\tau$, the function
${}_2F_1(\frac1{12},\frac5{12};1;\hat J(\tau))$ can be shown to be a
weight-$1$ modular form for~$\Gamma$, and to equal $\hat J^{-1/2}(\hat
J-1)^{-1/4}(d\hat J/d\tau)^{1/2}$.

The entire two-dimensional space of solutions of~(\ref{eq:lowest}), viewed
as functions of~$\tau$, is ${}_2F_1(\frac1{12},\frac5{12};1;\hat J(\tau))
\,(\mathbb{C}\tau +\mathbb{C})$ \,\cite{Stiller88}.  By comparison, the
space of solutions of~$D_\tau^2\tilde u=0$ on the $\tau$-plane
is~$\mathbb{C}\tau+\mathbb{C}$.  Lifting $L_{\frac1{12},\frac5{12};1}\,
u=0$ from $X(1)$ to~$\mathfrak{H}$ along the infinite-sheeted covering map
$\tau\mapsto \hat J$ will yield $D_\tau^2\tilde u=0$ only if the lifting is
`weak': if the dependent variable is altered according to $\tilde
u={}_2F_1(\frac1{12},\frac5{12};1;\hat J(\tau))^{-1}\cdot u$.

Proposition~\ref{prop:0} can be generalized from
$X(1)=\Gamma\setminus\mathfrak{H}^*$ to any
$X_1=\Gamma_1\setminus\mathfrak{H}^*$.  There are two cases to be
distinguished: when this algebraic curve is of genus zero, with function
field $\mathbb{K}(X_1)=\mathbb{C}(x)$ where $x$~is any Hauptmodul, and when
it is of positive genus, with function field
$\mathbb{K}(X_1)=\mathbb{C}(x,y)$, where $x,y$~are related by some
polynomial equation $\Phi(x,y)=0$ over~$\mathbb{C}$.

Before stating the generalization, we recall the definition of
characteristic exponents of a second-order differential operator
$L=D_x^2+\mathcal{A}\cdot D_x+\mathcal{B}$ on~$X_1$, where
$\mathcal{A},\mathcal{B}\in\mathbb{K}(X_1)$.  Such an operator is said to
be Fuchsian if all its singular points are regular, i.e., if it has two
characteristic exponents $\alpha_{i,1},\alpha_{i,2}\in\mathbb{C}$ (which
may be the same) at each singular point~$\mathfrak{s}_i\in X_1$.  If
$\alpha_{i,1}-\alpha_{i,2}\not\in\mathbb{Z}$, this means $L u=0$ has
local solutions $ u_{i,j}$, $j=1,2$, at~$\mathfrak{s}_i$ of the form
$t^{\alpha_{i,j}}$~times an invertible function of~$t$, where $t$~is a
local uniformizing parameter (if $\alpha_{i,1}-\alpha_{i,2}\in\mathbb{Z}$,
one solution may be logarithmic).  For example, at each singular point with
a zero exponent there is an {\em analytic\/} local solution, unique up~to
normalization.  Each ordinary (i.e., non-singular) point has
exponents~$0,1$.

\begin{theorem}
\label{thm:ford1}
In a neighborhood of any point on~$X_1=\Gamma_1\setminus\mathfrak{H}^*$,
any branch of the function~$\tau$ equals the ratio~$ u_1/ u_2$ of two
independent local solutions of a second-order Fuchsian differential
equation $L u:=(D_x^2+\mathcal{A}\cdot D_x+\mathcal{B}) u=0$, where
$\mathcal{A},\mathcal{B}\in\mathbb{K}(X_1)$.  Moreover, any such ratio is
of the form $(a\tau+b)/(c\tau+d)$ with $ad\neq bc$.  One can choose~$L$ so
that its singular points are the fixed points of\/~$\Gamma_1$, with the
difference of characteristic exponents equalling\/ $1/k$ at each fixed point
of order~$k$, and zero at each parabolic fixed point\/ {\rm(}i.e.,
cusp\/{\rm)}.
\end{theorem}

This is a special case of a theorem dealing with Fuchsian automorphic
functions of the first kind~\cite[\S\,44, Thm.~15]{Ford51}.  The
coefficients $\mathcal{A},\mathcal{B}$ of~$L$ can be taken to be
$0,\frac12\{\tau,x\}$, where $\{\cdot,\cdot\}$ is the Schwarzian
derivative.  In this case the space of solutions of $L u=0$, regarded as
functions of~$\tau$, is $(dx/d\tau)^{1/2}(\mathbb{C}\tau+\mathbb{C})$.
However, $L$~is not unique.  Any substitution of the type $ u\mapsto
f^\alpha u$, where $f\in\mathbb{K}(X_1)$ and~$\alpha\in\mathbb{C}$, will
produce an operator with transformed coefficients
$\hat{\mathcal{A}},\hat{\mathcal{B}}\in\mathbb{K}(X_1)$, but the same
solution ratios.  Similarly, the Liouvillian substitution
$ u\mapsto u\exp\left[-\frac12\smallint_{x_0}^x\mathcal{A}\,dx\right]$
will transform~$\mathcal{A},\mathcal{B}$, where $\mathcal{A}$~may be
nonzero,
to~$0,\mathcal{B}-\frac12D_x\mathcal{A}-\frac14{\mathcal{A}}^2$.
These substitutions will leave exponent {\em differences\/} invariant,
though they may shift exponents.

\begin{theorem}
\label{thm:ford2}
Let $L_i=D_x^2+\mathcal{A}_i\cdot D_x+\mathcal{B}_i$, $i=1,2$, where
$\mathcal{A}_i,\mathcal{B}_i\in\mathbb{K}(X_1)$, be Fuchsian operators
on~$X_1=\Gamma_1\setminus\mathfrak{H}^*$ with the property that any ratio
of independent solutions of $L_1 u=0$ or of~$L_2 u=0$ is of the
form~$(a\tau+b)/(c\tau+d)$.  Suppose that the singular points of each are
the fixed points of\/~$\Gamma_1$, with the difference of exponents
equalling\/ $1/k$ at each elliptic fixed point of order~$k$ and zero at
each cusp.  {\rm(}i\,{\rm)}~If $\mathcal{A}_1=\mathcal{A}_2=0$, then
$L_1,L_2$ are equal.  {\rm(}ii\,{\rm)}~If $X_1$~is of genus zero and
$L_1,L_2$ have the same\/ {\em exponents} {\rm(}not merely exponent
differences\/{\rm)}, then they are equal.
\end{theorem}

Part~(i) is a special case of a theorem dealing with Fuchsian automorphic
functions \cite[\S\,111, Thm.~7; and~\S\,115]{Ford51}.  The uniqueness
of~$L$ is deduced from the requirement that at each point, the inverse of
any solution ratio $x\mapsto u_1/ u_2$, such as $x$~as a function
of~$\tau\in\mathfrak{H}$, must be single-valued.  Part~(ii) is proved by
applying the above Liouvillian substitution to~$L_1,L_2$, and then invoking
part~(i).

To place this theorem in context, recall that any Fuchsian operator~$L$
on~$X_1$ of the form $D_x^2+\mathcal{A}\cdot D_x+\mathcal{B}$ determines a
flat (i.e., integrable) holomorphic connection on a trivial $2$-plane
bundle over the punctured curve
$X_1^0:=X_1\setminus\{\mathfrak{s}_1,\dots,\mathfrak{s}_n\}$, where
$\mathfrak{s}_1,\dots,\mathfrak{s}_n$ are the singular points of~$L$.  The
connection comes from analytically continuing any two independent local
solutions~$ u_1, u_2$ along paths in~$X_1^0$, producing an element of~${\it
GL}(2,\mathbb{C})$ for each path.  One may optionally quotient~out the
center~$\mathbb{C}^\times\!$ of~${\it GL}(2,\mathbb{C})$ to obtain a
projective action: an element of the M\"obius group ${\it
PGL}(2,\mathbb{C})$, acting on the ratio $ u_1/
u_2\in\mathbb{P}^1(\mathbb{C})$.  Such actions constitute a flat
holomorphic connection on a trivial $\mathbb{P}^1(\mathbb{C})$-bundle
over~$X_1^0$.

In~general, $L$~and the consequent flat connection on the trivial $2$-plane
bundle over~$X_1^0$ are not uniquely determined by the $2n-1$ independent
exponents.  In classical language, $L$~and the flat connection are
parametrized by the exponents together with certain (complex) {\em
accessory parameters\/}.  The projectivized flat connection on the trivial
$\mathbb{P}^1(\mathbb{C})$-bundle over~$X_1^0$ will be parametrized by the
$n$~exponent {\em differences\/}, together with certain of the accessory
parameters.  If $n\ge3$ and $X_1$~is of genus~$g$, the projectivized
connection will depend on $n-3+3g$ `projective' accessory parameters.  If
the innocuous normalization $\mathcal{A}=0$ is imposed, then the remaining
coefficient~$\mathcal{B}$ will be naturally parametrized by the exponent
differences and these accessory parameters.  Part~(i) of
Theorem~\ref{thm:ford2} is really a statement that if $x\mapsto u_1/ u_2$
is to be the inverse of a single-valued function at each point, then the
projective accessory parameters, and hence the flat connection on the
trivial $\mathbb{P}^1(\mathbb{C})$-bundle, are uniquely determined by the
exponent differences.

In~the absence of any imposed normalization, the
pair~$\mathcal{A},\mathcal{B}$, and hence $L$~itself and the flat
connection on the trivial $2$-plane bundle, will be parametrized by the
$2n-1$ independent exponents, the $n-3+3g$ projective accessory parameters,
and $g$~additional `affine' accessory parameters.  The presence of these
parameters when $g>0$ is the reason for the restriction to $g=0$ in
part~(ii) of Theorem~\ref{thm:ford2}.  Their values could differ between
$L_1$ and~$L_2$, even if those two differential operators have the same
exponents and projective accessory parameters.

The following is an explanation of how the accessory parameters, both
projective and affine, appear in any operator~$L$ of the form
$D_x^2+\mathcal{A}\cdot D_x+\mathcal{B}$.  If $X_1$~is of genus zero,
suppose that one of the $n\ge3$~singular points is~$x=\infty$, and that
each finite singular point has one exponent equal to zero.  (For instance,
the operator $L_{\frac{1}{12},\frac{5}{12};1}$ of
${}_2\mathcal{E}_1(\frac1{12},\frac5{12};1)$ has these properties.)  Then
$L$~will have the normal form~\cite{Poole36}
\begin{equation}
\label{eq:Poole}
\frac{d^2}{dx^2} +
\left[\sum_{i=1}^{n-1}\frac{1-\rho_i}{x-a_i}\right]\cdot\frac{d}{dx}
+\left[\frac{\Pi_{n-3}(x)}{\prod_{i=1}^{n-1}(x-a_i)}\right].
\end{equation}
Here $\{a_i\}_{i=1}^{n-1}$ are the finite singular points, with
exponents~$0,\rho_i$, and $\Pi_{n-3}(x)$ is a degree-$(n-3)$ polynomial.
Its leading coefficient determines the exponents at~$x=\infty$, and its
$n-3$~trailing coefficients are the (projective) accessory parameters.  If
on the other hand $X_1$~is of positive genus, suppose for simplicity it is
elliptic ($g=1$) or hyperelliptic ($g>1$), with $\Phi(x,y)=y^2-{\sf
P}_{2g+2}(x)$, where ${\sf P}_{2g+2}(x)$ is some polynomial of
degree~$2g+2$ with simple roots.  Then the generalization
of~(\ref{eq:Poole}) is straightforward.  If the $n$~singular points include
the infinite points $(x,y)=(\infty,\pm\infty)$ and $n-2$ finite points
$(x,y)=(x_i,y_i)$, the exponents of the latter being~$0,\rho_i$, then
$L$~will necessarily be proportional to
\begin{multline}
\label{eq:Pooleplus}
\left(y\frac{d}{dx}\right)^2 +
\left[\sum_{i=0}^{n-2}(1-\rho_i)\left(\frac12\cdot\frac{y+y_i}{x-x_i}\right) + \Pi_g(x)\right]\cdot\left(y\frac{d}{dx}\right)
\\{}+\left[
\sum_{i=1}^{n-2}{\mathfrak{b}}_i\left(\frac12\cdot\frac{y+y_i}{x-x_i}\right)
+\Pi_{2g}(x) +y\,\Pi_{g-1}(x)
\right],
\end{multline}
where the operator $y\,d/dx$ is the fundamental derivation
on~$\mathbb{C}(x,y)$, and the function $(1/2)(y+y_i)/(x-x_i)$ has simple
poles at $(x_i,y_i)$ and~$(\infty,\pm\infty)$.  The three polynomials
$\Pi_g(x)$, $\Pi_{2g}(x)$, $\Pi_{g-1}(x)$ have respective degrees
$g,2g,g-1$.  Their leading coefficients determine the exponents of the two
infinite points.  The $g$~trailing coefficients of~$\Pi_g(x)$ are the
affine accessory parameters, and the $n-2$ coefficients
$\{{\mathfrak{b}}_i\}_{i=1}^{n-2}$, the $2g$~trailing coefficients
of~$\Pi_{2g}(x)$, and the $g-1$~trailing coefficients of~$\Pi_{g-1}(x)$,
together make~up the $n-3+3g$ projective accessory parameters.

\subsection{Weak liftings}

We now specialize to the case when the Fuchsian group~$\Gamma_1$ is a
subgroup of the modular group~$\Gamma$, and in~particular to the case
$\Gamma_1=\Gamma_0(N)$.  By Theorem~\ref{thm:ford1}, it is possible to
represent the coordinate~$\tau$ of~$\mathfrak{H}$ as a ratio of two
solutions of some Fuchsian differential equation
on~$X_0(N)=\Gamma_0(N)\setminus\mathfrak{H}^*$.  A~differential equation
with this property can be derived as the lifting, or any weak lifting,
of~${}_2\mathcal{E}_1(\frac1{12},\frac5{12};1)$ from $X(1)$ to~$X_0(N)$.
Under some circumstances, as will be explained in the next section, it is
possible to derive a similar equation on~$X_0(N^2)$ by either of two
further liftings, the equivalence between which will yield a hypergeometric
identity.  Any two such equations are identical only if their exponents are
the same; so we now consider the effects of liftings on exponents.

Suppose $L=D_x^2+\mathcal{A}\cdot D_x+\mathcal{B}$ is a differential
operator on any algebraic curve $X_1=\Gamma_1\setminus\mathfrak{H}^*$
over~$\mathbb{C}$ with derivation~$D_x$ and function
field~$\mathbb{K}/\mathbb{C}$.  Let $\xi:\tilde X_1\to X_1$ be the rational
map of curves arising from a subgroup relation~$\tilde\Gamma_1\le\Gamma_1$.
Here $\tilde X_1$~will have its own function
field~$\tilde{\mathbb{K}}/\mathbb{C}$ with derivation~$D_{\tilde x}$.  The
lifting of~$L$ to~$\tilde X_1$ will be an operator $\tilde L=D_{\tilde x}^2
+ \tilde{\mathcal{A}}\cdot D_{\tilde x}+\tilde{\mathcal{B}}$ with
$\tilde{\mathcal{A}},\tilde{\mathcal{B}}\in\tilde{\mathbb{K}}$, satisfying
the condition that $Lu=0$, $\tilde L\tilde u=0$ locally have independent
solution pairs $u_1,u_2$ and~$\tilde u_1,\tilde u_2$ such that $\tilde
u_i=u_i\circ\xi$.  Informally, $\tilde L$~is obtained from~$L$ by
performing a change of (independent) variable, and also left-multiplying by
an element of~$\tilde{\mathbb{K}}^\times\!$ if needed to preserve monicity.
Each local solution~$\tilde u(\cdot)$ of $\tilde L\tilde u=0$ on~$\tilde
X_1$ will be of the form $u(\xi(\cdot))$, where $u$~is some local solution
of~$Lu=0$.

A differential operator~$\tilde M$ of the same form as~$\tilde L$ is said
to be a {\em weak\/} lifting of~$L$ if there are {\em ratios\/}
$\sigma,\tilde\sigma$ of independent solutions of $Lu=0$, $\tilde M\tilde
u=0$ respectively, such that $\tilde\sigma=\sigma\circ\xi$.  In~a
projective context, weak liftings are clearly the more appropriate concept.
Informally, a weak lifting of the differential equation $Lu=0$ to~$\tilde
X_1$ may incorporate a linear change of the {\em dependent\/} variable.
For example, it may be multiplied by any~$\tilde
f\in\tilde{\mathbb{K}}^\times\!$.  In this case any local solution~$\tilde
u(\cdot)$ of $\tilde M\tilde u=0$ on~$\tilde X_1$ will be of the form
$\tilde f(\cdot)\,u(\xi(\cdot))$, where $u$~is a local solution of~$Lu=0$.

From this point, we shall consider only lifting prefactors that are
elements of~$\tilde{\mathbb{K}}_e^\times\!
\supset\tilde{\mathbb{K}}^\times\!$, the collection of (algebraic,
multivalued) functions on~$\tilde X_1$ of the form $\tilde f^{\alpha}$ with
$\tilde f\in{\tilde {\mathbb{K}}}^\times\!$ and~$\alpha\in\mathbb{Q}$.  The
`extension' $\tilde{\mathbb{K}}_e^\times\!$ is not a field, due~to the
absence of closure under addition.  Up~to scalar multiplication, any
$\tilde f\in\tilde{\mathbb{K}}^\times\!$ may be identified with its divisor
$\sum_in_i\tilde{\mathfrak{p}}_i$, an element of the free
$\mathbb{Z}$-module on~$\tilde X_1$.  In~consequence, any~$\tilde
f_e\in\tilde{\mathbb{K}}_e^\times\!$ has an associated `generalized
divisor' of the form $\sum_i c_i\tilde{\mathfrak{p}}_i$, an element of the
free $\mathbb{Q}$-module on~$\tilde X_1$.  Since $(\tilde f^\alpha)\circ
D_{\tilde x}\circ (\tilde f^\alpha)^{-1}=D_{\tilde x}-\alpha(D_{\tilde
x}\tilde f)/\tilde f$, multiplying the dependent variable by any $\tilde
f_e\in\tilde{\mathbb K}_e^\times\!$ will yield a weak lifting~$\tilde M$
that has the same general form as~$\tilde L$, with coefficient
functions~$\tilde{\mathcal{A}},\tilde{\mathcal{B}}$ that are elements
of~$\tilde{\mathbb{K}}$.

Lifting a Fuchsian operator $L=D_x^2+\mathcal{A}\cdot D_x+\mathcal{B}$
on~$X_1$ to~$\tilde X_1$ transforms exponents in a straightforward way.
Suppose $\tilde{\mathfrak{p}}\in\tilde X_1$ is a critical point
of~$\xi:\tilde X_1\to X_1$ with corresponding critical
value~$\mathfrak{p}\in X_1$, i.e., suppose $\xi(\tilde t)$ equals
$t^n$~times an invertible function of~$t$, where $\tilde t,t$ are local
uniformizing parameters near~$\tilde{\mathfrak{p}},\mathfrak{p}$
respectively, and $n>1$~is the ramification index.  Then the exponents of
the lifting~$\tilde L$ at~$\tilde{\mathfrak{p}}\in\tilde X_1$ will be
$n$~times those of~$L$ at~$\mathfrak{p}\in X_1$.  (This statement extends
to the case when $\tilde{\mathfrak{p}}$~is not a critical point of~$\xi$,
and~$n=1$.)  If~the lifting prefactor~$\tilde
f_e\in\tilde{\mathbb{K}}_e^\times\!$ has generalized divisor
$\sum_ic_i\tilde{\mathfrak{p}}_i$, where $c_i\in\mathbb{Q}$
and~$\tilde{\mathfrak{p}}_i\in\tilde X_1$, then the exponents of the weak
lifting~$\tilde M$ at each point~$\tilde{\mathfrak{p}}_i$ will be shifted
by~$c_i$, relative to those of the lifting~$\tilde L$.

Most previous work on Fuchsian differential equations on modular curves has
adhered to a convention taken from conformal mapping, according to which
the equation should be of the form $(D_x^2+{\mathcal{B}}) u=0$, i.e.,
should be formally self-adjoint, with~${\mathcal A}=0$~\cite{Harnad2000}.
This typically forces each exponent to be nonzero.  To~derive
hypergeometric identities, it~is better to adopt an asymmetric convention
informally introduced in~(\ref{eq:Poole})--(\ref{eq:Pooleplus}): every
singular point but one should have a zero exponent.

\begin{definition}
A Fuchsian operator $L=D_x^2+\mathcal{A}\cdot D_x + \mathcal{B}$ on an
algebraic curve~$X_1$ over~$\mathbb{C}$ is said to be in normal form
relative to a specified point~$\varinfty\in X_1$ (typically, a singular
point) if each singular point of~$L$ not equal to~$\varinfty$ has a zero
exponent.
\end{definition}

\noindent
Requiring a weak lifting to be in normal form comes close to specifying it
uniquely.  This will not be explored further here, since for the purpose of
studying weak liftings from one modular curve to another, the following
lemma will suffice.

\begin{lemma}
\label{lem:L}
Suppose a Fuchsian operator $L=D_x^2+\mathcal{A}\cdot D_x + \mathcal{B}$
on~$X_1$ is in normal form relative to a specified point~$\varinfty\in
X_1$, with $\alpha_\infty\!$ denoting one of the exponents of~$L$
at\/~$\varinfty$.  Let $\xi:\tilde X_1\to X_1$ be a rational map of
algebraic curves over\/~$\mathbb C$.  If $\tilde f\in\tilde{\mathbb{K}}$ is
a function with divisor equal to\/
$\sum_{\tilde{\mathfrak{p}}\in\xi^{-1}(\varinfty)}\left[(\tilde{\mathfrak{p}})-(\widetilde\varinfty)\right]$,
for some\/ $\widetilde\varinfty\in\xi^{-1}(\varinfty)$, then the weak
lifting $\tilde M= D_{\tilde x}^2+\tilde{\mathcal{A}}\cdot D_{\tilde
x}+\tilde{\mathcal{B}}$ of~$L$ to~$\tilde X_1$ induced by $\tilde
u(\cdot)=\tilde f(\cdot)^{-\alpha_\infty}\,u(\xi(\cdot))$ will be in normal
form relative to\/~$\widetilde\varinfty$.

Moreover, for any singular point $\tilde O\notin \xi^{-1}(\varinfty)$
on~$\tilde X_1$, the unique local solution~$\tilde u_0$ of~$\tilde M\tilde
u=0$ which is analytic and equal to unity at~$\tilde O$ equals $\tilde
f(\cdot)^{-\alpha_\infty}u_0(\xi(\cdot))$, where $u_0$~is the corresponding
analytic local solution of~$Lu=0$ at~$O:=\xi(\tilde O)$.
\end{lemma}
\begin{proof}
The first statement follows from a straightforward computation of the
lifted exponents; and the second from the preservation of analyticity of
local solutions at each singular point other than~$\underline{\infty}$,
under lifting.
\end{proof}

\begin{lemma}
\label{lem:L0}
Let $\xi:\mathbb{P}^1(\mathbb{C})_{\tilde x}\to\mathbb{P}^1(\mathbb{C})_x$
be a rational map from the $\tilde x$-sphere to the $x$-sphere, and suppose
the operator $L=D_x^2+\mathcal{A}(x)D_x+\mathcal{B}(x)$ on the $x$-sphere
is in normal form relative to the point~$x=\infty$.  Let $\alpha_\infty$
denote one of the two exponents of~$Lu=0$ at~$x=\infty$.  If $\xi(\tilde
x)=P(\tilde x)/Q(\tilde x)$ with $P,Q$ having no nontrivial factor
in~common, then the weak lifting~$\tilde M$ of~$L$ to the $\tilde x$-sphere
induced by $\tilde u(\cdot)=Q(\cdot)^{-\alpha_\infty}u(\xi(\cdot))$ will be
in normal form relative to the point~$\tilde x=\infty$.

Moreover, for any singular point $\tilde x=\tilde O$ with $\xi(\tilde
O)\neq\infty$, the unique local solution~$\tilde u_0$ of~$\tilde M\tilde
u=0$ analytic and equal to unity at~$\tilde x=\tilde O$ equals
$Q(\cdot)^{-\alpha_\infty}\,u_0(\xi(\cdot))$, where $u_0$~is the
corresponding analytic local solution of~$Lu=0$ at~$O:=\xi(\tilde O)$.
\end{lemma}
\begin{proof}
This is a specialization of Lemma~\ref{lem:L} to the case of zero genus.
\end{proof}

As the appendix summarizes, if $g(X_0(N))=0$ then $X_0(N)$~is coordinatized
by a Hauptmodul~$x_N$ with divisor $({\rm i}\infty)-(0)$, i.e., a~univalent
function~$x_N$ with a simple zero (resp.~pole) at the cusp~$\tau={\rm
i}\infty$ (resp.~the cusp~$\tau=0$).  The covering map $\pi_N:X_0(N)\to
X(1)$ is given by a rational function $j=j(x_N)=P_N(x_N)/Q_N(x_N)$,
a~quotient of monic polynomials satisfying $\deg P_N=\psi(N)$, $\deg
P_N-\deg Q_N=N$, and $x_N|Q_N(x_N)$.  So $\hat J=\hat
J(x_N)=1728\,Q_N(x_N)/P_N(x_N)$.  The cusps on the $x_N$-sphere $X_0(N)$
comprise the zeroes of~$Q_N$, including~$x_N=0$; and also~$x_N=\infty$.

\begin{definition}
\label{def:hNdef}
If $X_0(N)$ is of genus zero, the fundamental analytic function~$h_N$,
which will play a major role in the sequel, is defined by
\begin{equation}
\label{eq:hNdef}
h_N(x_N):=\left[P_N(x_N)/P_N(0)\right]^{-1/12}\!{}_2F_1\left(\tfrac1{12},\tfrac5{12};\,1;\,
\frac{1728\,Q_N(x_N)}{P_N(x_N)}\right)
\end{equation}
in a neighborhood of the distinguished cusp $x_N=0$ of~$X_0(N)$, i.e.,
of~$\tau={\rm i}\infty$.
\end{definition}

\begin{proposition}
\label{prop:2}
If $X_0(N)$ is of genus zero, there is a normal-form weak lifting $\tilde
M_N\tilde u=0$ of ${}_2\mathcal{E}_1(\frac1{12},\frac5{12};1)$ from the
$\hat J$-sphere~$X(1)$ to~$X_0(N)$, along~$\pi_N$, which has a total of
$(\sigma_\infty+\varepsilon_{\rm i}+\varepsilon_\rho)(N)$ singular points
on~$X_0(N)$, classified thus:
\begin{enumerate}
\item one singular point with characteristic
exponents\/~$\frac1{12}\psi(N),\frac1{12}\psi(N)$, namely the cusp
${x_N=\infty}$ {\rm(}i.e.,~the distinguished cusp~$\tau=0${\rm)};
\item $\sigma_\infty(N)-1$~singular points with exponents\/~$0,0$, namely
the remaining cusps, including $x_N=0$ {\rm(}i.e., the distinguished
cusp~$\tau={\rm i}\infty${\rm)};
\item $\epsilon_{\rm i}(N)$~singular points with exponents\/~$0,\frac12$,
namely the order-$2$ elliptic fixed points;
\item $\epsilon_{\rho}(N)$~singular points with exponents\/~$0,\frac13$,
namely the order-$3$ elliptic fixed points.
\end{enumerate}
This weak lifting is obtained from the lifting prefactor $\tilde
f_e(x_N)=P_N(x_N)^{-1/12}$.

Any ratio of independent local solutions of~$\tilde M_N\tilde u=0$ equals
$(a\tau+b)/(c\tau+d)$ for some $a,b,c,d$ with $ad\neq bc$.  The unique
local solution of~$\tilde M_N\tilde u=0$ analytic at the distinguished cusp
$x_N=0$ {\rm(}i.e., at~$\tau={\rm i}\infty${\rm)} and equalling\/~$1$ there
will be~$h_N$.
\end{proposition}
\begin{remarkaftertheorem}
The differential equation $\tilde M_N\tilde u=0$ satisfies the conditions
of Theorem~\ref{thm:ford1}. By Theorem~\ref{thm:ford2}(ii), it is uniquely
characterized by the given list of singular points and exponents, and the
fact that any ratio of independent local solutions is of the form
$(a\tau+b)/(c\tau+d)$.  Its two-dimensional space of solutions, viewed as
functions of~$\tau$, is $h_N(x_N(\tau))\,(\mathbb{C}\tau +\mathbb{C})$.
\end{remarkaftertheorem}
\begin{proof}[Proof of Proposition\/~{\rm\ref{prop:2}}]
This is an application of Lemma~\ref{lem:L0}, with
$\alpha_\infty=\frac1{12}$.  The given exponents are computed thus.
Suppose the prefactor were absent, i.e., suppose $\tilde
u(\cdot)=u(\xi(\cdot))$.  Then above $\hat J=0$, each cusp on~$X_0(N)$
would have exponents~$0,0$, irrespective of its multiplicity.  A~point on
the fibre above~$\hat J=1$ would have exponents~$0,\frac12$ if it has unit
multiplicity, i.e., if it is an order-$2$ elliptic point; and
exponents~$0,1$ otherwise, indicating it would be an ordinary
(non-singular) point.  The simple roots of~$P_N(x_N)$ are the
$\varepsilon_\rho(N)$ order-$3$ elliptic points on the $x_N$-sphere, and
all other roots are triple.  For a lifting, the exponents at these two
sorts of point would be $\frac1{12},\frac5{12}$ and~$\frac14,\frac54$.
Including the prefactor $P_N(x_N)^{-1/12}$ shifts them to~$0,\frac13$
and~$0,1$; so the latter will no~longer be singular points.  It~also shifts
the exponents at the cusp~$x_N=\infty$ from~$0,0$ to~$\frac1{12}\deg
P_N,\frac1{12}\deg P_N$.
\end{proof}

The weak-lifted differential equation $\tilde M_N\tilde u=0$ of the
proposition, which is of the form $\left[D_{x_N}^2+\tilde{\mathcal{A}}\cdot
D_{x_N}+\tilde{\mathcal{B}}\right]\tilde u=0$ for certain
$\tilde{\mathcal{A}},\tilde{\mathcal{B}}\in\mathbb{Q}(x_N)$, can readily be
derived from ${}_2\mathcal{E}_1(\frac1{12},\frac5{12};1)$ by changing
variables.  Like ${}_2\mathcal{E}_1(\frac1{12},\frac5{12};1)$, it is always
based on an operator of the normal form~(\ref{eq:Poole}).  The
curve~$X_0(N)$ is of genus zero only if
$N=2,3,4,5,6,7,8,9,10,12,13,16,18,25$, and the number of singular points,
namely $m(N):=\left(\sigma_\infty+\varepsilon_{\rm
i}+\varepsilon_\rho\right)(N)$, is $3,3,3,4,4,4,4,4,6,6,6,6,8,8$,
respectively.  So when $N=2,3,4$, the equation $\tilde M_N\tilde u=0$ is of
hypergeometric type on the $x_N$-sphere.  (Two of its singular points are
at the cusps $x_N=0,\infty$, so a linear scaling of~$x_N$ will reduce it to
the Gau\ss\ form.)  When~$N=5,6,7,8,9$, it~is of Heun type, and its
solutions, including~$h_N$, may be expressed in~terms of the local Heun
function~$\Hl$~\cite{Ronveaux95}.  In~general, the coefficients
$\bigl\{c^{(N)}_n\bigr\}_{n=0}^\infty$ of the power series expansion
$h_N(x_N)=\sum_{n=0}^\infty c^{(N)}_nx_N^n$ will satisfy an
$\left[m(N)-1\right]$-term recurrence relation.

\begin{proposition}
The fundamental analytic function $h_N$ defined by\/~{\rm(\ref{eq:hNdef})}
when $X_0(N)$~is of genus zero, on a neighborhood of the point $x_N=0$
on~$X_0(N)$, extends by continuation to a weight-$1$ modular form
for\/~$\Gamma_0(N)$ on~$\mathfrak{H}^*\ni\tau$, with some multiplier system.
One may write
\begin{equation}
\label{eq:althNdef}
h_N(x_N(\tau)) = P_N(0)^{1/12} Q_N(x_N(\tau))^{-1/12}\,\eta^2(\tau),
\end{equation}
where $\eta(\cdot)$~is the Dedekind eta function.  This modular form is
regular and non-vanishing at each cusp of\/~$\Gamma_0(N)$
in\/~$\mathbb{P}^1(\mathbb{Q})\ni\tau$ other than those in the equivalence
class\/~$\left[\frac11\right]_N\ni0$, at each of which its order
is~$\psi(N)/12N$.
\end{proposition}
\begin{proof}
Stiller~\cite{Stiller88} shows that in a neighborhood of~$\tau={\rm
i}\infty$, where $\hat J=\hat J(\tau)$ equals zero, the analytic solution
${}_2F_1(\frac1{12},\frac5{12};1;\hat J)$ of
${}_2\mathcal{E}_1(\frac1{12},\frac5{12};1)$, regarded like~$\hat J$ as a
function of~$\tau$, equals $12^{1/4}\eta^2(\tau)\hat J^{-1/12}$, a
weight-$1$ modular form for~$\Gamma$.  Equivalently,
\begin{equation}
\label{eq:stiller}
\eta^2(\tau) = 12^{-1/4} \hat J^{1/12}\,
{}_2F_1(\tfrac1{12},\tfrac5{12};1;\hat J),
\end{equation}
with the root ${\hat J}^{1/12}$ taken positive when $\hat J>0$.  (The
identity (\ref{eq:stiller})~can be traced back to Dedekind, who at one
point {\em defined\/} $\eta^2(\tau)$ as an expression of hypergeometric
type equal to the right-hand side~\cite[p.~137]{Chandrasekharan85}.)
Combining this fact with $\hat J(x_N)=1728\,Q_N(x_N)/P(x_N)$
and~(\ref{eq:hNdef}) yields~(\ref{eq:althNdef}).

$\eta(\tau),x_N(\tau)$ are of weight $1/2,0$, and the $Q_N$~factor
in~(\ref{eq:althNdef}) is zero only at cusps of~$\mathfrak{H}^*$, so the
right-hand side of~(\ref{eq:althNdef}) is holomorphic on~$\mathfrak{H}$ and
is a weight-$1$ modular form.  The roots of the polynomial~$Q_N(x_N)$ are
bijective with the cusp equivalence classes $\left[\frac{a}d\right]_N$
of~$X_0(N)$ other than~$\left[\frac11\right]_N$, at~which $x_N=\infty$.
Any root of~$Q_N$ has multiplicity equal to the cusp width~$e_{d,N}$, i.e.,
the multiplicity with which $\left[\frac{a}d\right]_N$~is mapped to~$X(1)$.
Also, the order of~$\eta(\tau)$ at any
cusp~$\tau\in\mathbb{P}^1(\mathbb{Q})$ equals~$1/24$.  So altogether, the
order of~$h_N(x_N(\tau))$ at any cusp
in~$\left[\frac{a}d\right]_N\neq\left[\frac11\right]_N$ will equal
$(-1/12)e_{d,N}/e_{d,N}+2\cdot\frac1{24}=0$.

The cusp equivalence class~$\left[\frac11\right]_N$, on which~$x_N=\infty$,
is best handled by referring to the original definition~(\ref{eq:hNdef})
of~$h_N$.  Since $\deg P_N=\psi(N)$, the order of~$h_N(x_N(\tau))$ at any
cusp in~$\left[\frac11\right]_N$ will be
$(1/12)\psi(N)/e_{1,N}=\psi(N)/12N$.
\end{proof}

\begin{table}
\caption[The Hauptmodul~$x_N$ and weight-$1$ modular form~$h_N$
for~$Gamma_0(N)$, as eta products.]{The Hauptmodul~$x_N$ and weight-$1$
modular form~$h_N$ for~$\Gamma_0(N)$, as eta products.  (For notation see
the appendix.)}  \hfill
\begin{tabular}{|c|l|l|}
\hline
$N$ & \hfil $x_N(\tau)$ & \hfil $h_N(x_N(\tau))$\\
\hline
$2$ & $2^{12}\cdot[2]^{24}/\,[1]^{24}$ & $[1]^4/\,[2]^2\vphantom{\left\{[1]^5/\,[5]\right\}^{1/2}}$\\
$3$ & $3^{6}\cdot[3]^{12}/\,[1]^{12}$ & $[1]^3/\,[3]\vphantom{\left\{[1]^5/\,[5]\right\}^{1/2}}$\\
$4$ & $2^{8}\cdot[4]^{8}/\,[1]^{8}$ & $[1]^4/\,[2]^2\vphantom{\left\{[1]^5/\,[5]\right\}^{1/2}}$\\
$5$ & $5^{3}\cdot[5]^{6}/\,[1]^{6}$ & $\left\{[1]^5/\,[5]\right\}^{1/2}$\\
$6$ & $2^{3} 3^{2}\cdot [2][6]^{5}/\,[1]^{5}[3]$ & $[1]^6[6]\,/\,[2]^3[3]^2\vphantom{\left\{[1]^5/\,[5]\right\}^{1/2}}$\\
$7$ & $7^{2}\cdot[7]^{4}/\,[1]^{4}$ & $\left\{[1]^7/\,[7]\right\}^{1/3}$\\
\hline
\end{tabular}
\hfill
\label{tab:1}
\end{table}

For each of $N=2,\dots,7$, an eta-product representation for the weight-$1$
form $h_N(x_N(\tau))$ can be computed from~(\ref{eq:althNdef}), and the
formula $j=P_N(x_N)/Q_N(x_N)$ and eta-product representation for
$x_N=x_N(\tau)$ given in the appendix.  These are listed in
Table~\ref{tab:1}.  The multiplier systems of the $h_N(x_N(\tau))$ are
nontrivial but are not difficult to compute.  For example, $h_6(x_6(\tau))$
is of quadratic Nebentypus: its transformation under any
$\left(\begin{smallmatrix}a&b\\c&d\end{smallmatrix}\right)\in\Gamma_0(6)$
includes a $(\pm1)$-valued factor equal to the Dirichlet character
$\chi(d):=\left(\frac{d}6\right)$, where $\left(\frac{\cdot}{\cdot}\right)$
is the Jacobi symbol.  (Cf.~\cite[\S\,2.2]{Verrill99}.)

The new interpretation of~$h_N$ makes contact with the theory of
differential equations satisfied by modular forms, originating with
Stiller~\cite{Stiller84}.  If~$\Gamma_1\le\Gamma$ is a congruence subgroup
of genus zero with Hauptmodul~$x$, then any weight-$k$ modular form~$f$
for~$\Gamma_1$ that has a power series expansion $f=\sum_{n=0}^\infty
c_nx^n$ and can be viewed locally (near~$x=0$) as a function of~$x$ will
necessarily satisfy an order-$(k+1)$ differential equation with respect
to~$x$, which can be constructed algorithmically.  Applied to any~$h_N$
viewed as a weight-$1$ modular form for~$\Gamma_0(N)$, the algorithm will
recover the weak-lifted second-order equation $\tilde{M}_N\tilde u=0$ of
Proposition~\ref{prop:2}.

We shall not attempt here to define an analogue of the modular form
$h_N(x_N(\tau))$ when $X_0(N)$~is of positive genus.  Although the curves
$X_0(36),X_0(49)$ are of genus~$1$, it will nonetheless prove possible
in~\S\S\,\ref{sec:keyresults} and~\ref{sec:examples} to derive
hypergeometric identities from the way in which they cover $X_0(6),X_0(7)$,
at the price of some awkwardness.

\section{Key Results}
\label{sec:keyresults}

It can now be explained how hypergeometric identities follow from the
covering maps $\phi_N:X_0(N^2)\to X_0(N)$ and~$\phi_N':X_0(N^2)\to
X_0(N)'$.  These coverings are induced by the subgroup relations
${\Gamma_0(N^2)<\Gamma_0(N)}$ and~$\Gamma_0(N^2)<\Gamma_0(N)'$.  Note that
the map~$\phi_N$ also has an obvious elliptic-curve interpretation: it~acts
on~$X_0(N^2)$, the space of isomorphism classes of $N^2$-isogenies between
elliptic curves $E,E'=E/C_{N^2}$, by replacing $C_{N^2}$ by its $C_N$
subgroup.  The map~$\phi_N'$ acts on the dual $N$-isogeny
between~$E',E=E'/C_{N^2}$ in the same way.

If $X_0(N)$ is of genus zero, and has a canonical Hauptmodul~$x_N$ with
divisor $({\rm i}\infty)-(0)$ as reviewed in the appendix, then the
fundamental weight-$1$ modular form~$h_N$ is given near the point~$x_N=0$,
i.e., the cusp~$\tau={\rm i}\infty$, by Definition~\ref{def:hNdef}.
Consider first the case when $X_0(N^2)$ too is of genus zero, with
Hauptmodul~$x_{N^2}$, and a corresponding weight-$1$ modular form~$h_{N^2}$
defined near $x_{N^2}=0$ (which is the same cusp~$\tau={\rm i}\infty$
on~$\mathfrak{H}^*$).

Both $x_N$ and~$x'_N$, the canonical Hauptmodul on~$X_0(N)'$, will be
rational functions of~$x_{N^2}$.  That~is, $x_N=R(x_{N^2})/S(x_{N^2})$
and~$x_N'=R'(x_{N^2})/S'(x_{N^2})$ for certain
$R,S,R',S'\in\mathbb{Q}[x_{N^2}]$.  (Primes do~not indicate derivatives.)
Both rational functions will have mapping degree equal to~$N$, the index
of~$\Gamma_0(N^2)$ in~$\Gamma_0(N)$ or~$\Gamma_0(N)'$.

By Lemma~\ref{lem:cusplifting}, the fibre of~$\phi_N$ above $x_N=\infty$,
i.e, above the cusp~$\left[\frac11\right]_N$ of~$X_0(N)$
containing~$\tau=0$, consists of the point $x_{N^2}=\infty$, i.e., the
cusp~$\left[\frac11\right]_{N^2}$ of~$X_0(N^2)$, with multiplicity~$N$.  So
$\deg R=N$ and~$\deg S=0$; one may take $S(x_{N^2}):=1$.  Also by
Lemma~\ref{lem:cusplifting}, the fibre of~$\phi_N'$ above~$x_N'=0$, i.e.,
above the cusp~$\left[\frac1N\right]_N/N$ of~$X_0(N)'$, consists of the
point~$x_{N^2}=0$, i.e., the cusp~$\left[\frac1{N^2}\right]_{N^2}$, with
multiplicity~$N$.  Moreover, the fibre of~$\phi_N'$ above~$x_N'=\infty$,
i.e., above~$\left[\frac11\right]_N/N$, contains the
point~$x_{N^2}=\infty$, i.e., the cusp~$\left[\frac11\right]_{N^2}$, with
multiplicity~$1$.  So $R'(x_{N^2})\propto x_{N^2}^N$ and $\deg S'=N-1$; one
may take~$R'(x_{N^2}):= x_{N^2}^N$.  It turns~out that the two remaining
polynomials~$R,S'$ are always monic and in~$\mathbb{Z}[x_{N^2}]$, but that
fact will not be needed.

\begin{theorem}
\label{thm:1}
If $X_0(N)$ and $X_0(N^2)$ are of genus zero\/ {\rm(}i.e., if
$N=2,3,4,5${\rm)}, then the hypergeometric identity
\begin{displaymath}
h_{N^2}(x_{N^2})=h_N\left(R(x_{N^2})\right)=
\left[S'(x_{N^2})/S'(0)\right]^{-\psi(N)/12}\,
h_N\left(\frac{x_{N^2}^N}{S'(x_{N^2})}\right)
\end{displaymath}
holds in a neighborhood of the point~$x_{N^2}=0$ on~$X_0(N^2)$, i.e., of
the cusp~$\tau={\rm i}\infty$.
\end{theorem}
\begin{remarkaftertheorem}
The first equality says that as functions of~$\tau\in\mathfrak{H}^*$, the
modular forms $h_{N^2}(x_{N^2}(\tau))$ and~$h_N(x_N(\tau))$ are equal.
This is confirmed in the $N=2$ case by Table~\ref{tab:1}, and is easily
seen to be true of the other genus-zero cases~$N=3,4,5$.
\end{remarkaftertheorem}
\begin{proof}
${}_2\mathcal{E}_1(\frac1{12},\frac5{12};1)$ is weak-lifted by~$\pi_N$
to~$\tilde M_N\tilde u=0$ and by~$\pi_{N^2}$ to~$\tilde M_{N^2}\tilde u=0$.
The unique (normalized) analytic local solutions at~$x_N=0$, resp.\
$x_{N^2}=0$, are $h_N$ and~$h_{N^2}$.  The expression $h_N(R(\cdot))$ is
the unique (normalized) analytic local solution at~$x_{N^2}=0$ of what may
be denoted $\tilde{\tilde M}_{N^2}\tilde u=0$, the {\em lifting\/}
(not~merely weak lifting) of~$\tilde M_N \tilde u=0$ from $X_0(N)$
to~$X_0(N^2)$, along~$\phi_N$.  Since $\pi_{N^2}=\pi_N\circ\phi_N$, the
equation $\tilde{\tilde M}_{N^2}\tilde u=0$ like $\tilde M_{N^2}\tilde u=0$
is a weak lifting of~${}_2\mathcal{E}_1(\frac1{12},\frac5{12};1)$
to~$X_0(N^2)$, along~$\pi_{N^2}$.  But $\tilde M_{N^2},\tilde{\tilde
M}_{N^2}$ have the same exponents: e.g.,
$N\left(\frac1{12}\psi(N),\frac1{12}\psi(N)\right)=\left(\frac1{12}\psi(N^2),\frac1{12}\psi(N^2)\right)$
at the cusp~$x_{N^2}=\infty$, and $0,0$~at each of the other cusps.  So by
Theorem~\ref{thm:ford2}(ii) they are the same, and $h_N(R(\cdot))$ must
equal~$h_{N^2}(\cdot)$.

To show the second expression also equals $h_{N^2}(x_{N^2})$, consider
$\tilde M'_N\tilde u=0$, the Fuchsian equation on the
$x_N'$-sphere~$X_0(N)'$ obtained by formally substituting $x_N'$ for~$x_N$
in~$\tilde M_N\tilde u=0$, and its lifting $\tilde{\tilde M}_{N^2}'\tilde
u=0$ to~$X_0(N^2)$, along~$\phi_N'$.  Any ratio of independent local
solutions of~$\tilde M_N\tilde u=0$ is of the form $(a\tau+b)/(c\tau+d)$.
Since $x_N'(\tau):=x_N(N\tau)$, the same is true of~$\tilde M_N'\tilde
u=0$.  (The coefficients~$a,c$ are multiplied by~$N$.)  By the definition
of a lifting, it is also true of~$\tilde{\tilde M}_{N^2}'\tilde u=0$; and
by Proposition~\ref{prop:2} it happens to be true of~$\tilde M_{N^2}\tilde
u=0$, as~well.  So $\tilde M_{N^2},\tilde{\tilde M}'_{N^2}$ are
projectively the same: they determine the same flat connection on the
$\mathbb{P}^1(\mathbb{C})$-bundle over~$X_0(N^2)$.  They necessarily have
the same exponent {\em differences\/}.

To make $\tilde{\tilde M}'_{N^2}$ identical to~$\tilde M_{N^2}$, it~may be
redefined as a {\em weak\/} lifting, incorporating a prefactor that shifts
the exponents at each singular point to those of~$\tilde M_{N^2}$.  (By
Theorem~\ref{thm:ford2}(ii), that will suffice.)  Only the exponents at the
cusps of~$X_0(N^2)$ need to be altered.  As~originally defined,
$\tilde{\tilde M}'_{N^2}$~had exponents
$n\left(\frac1{12}\psi(N),\frac1{12}\psi(N)\right)$ at each cusp in the
fibre of~$\phi_N'$ above~$x_N'=\infty$, where $n$~is the multiplicity with
which the cusp appears.  Letting the lifting prefactor be
$S'(x_{N^2})^{-\psi(N)/12}$, which is a member of the collection of
algebraic functions on the upper curve~$X_0(N^2)$
denoted~$\tilde{\mathbb{K}}_e^\times\!$ in the last section, will shift the
exponents at all cusps other than~$x_{N^2}=\infty$ to~$0,0$ as~desired.  So
with this choice, $\tilde{\tilde M}'_{N^2}$ will equal~$\tilde M_{N^2}$,
and the second expression of the theorem, which is the unique (normalized)
analytic local solution of $\tilde{\tilde M}'_{N^2}\tilde u=0$
at~$x_{N^2}=0$, will equal $h_{N^2}(x_{N^2})$, the corresponding local
solution of~$\tilde M_{N^2}\tilde u=0$.
\end{proof}

If $X_0(N^2)$ unlike~$X_0(N)$ is of positive genus, then it will have
no~Hauptmodul~$x_{N^2}$.  In this case Theorem~\ref{thm:1} must be replaced
by Theorem~\ref{thm:2}, which actually subsumes~it, though it~is less
explicit and must be supplemented by the subsequent proposition.  But what
it says about the case $g(X_0(N^2))>0$ is significantly weaker than
Theorem~\ref{thm:1}.  If the genus is positive, Theorem~\ref{thm:ford2}(ii)
is not available, and it is consequently difficult to rule~out the
possibility that the two sides of the hypergeometric `identity' may in~fact
differ, due~to their being solutions of Fuchsian differential equations
with different values for their affine accessory parameters.  In~the cases
$N=6,7$ treated in the next section, this will fortunately not be an issue,
since the two differential equations on~$X_0(N^2)$ can be computed
explicitly and shown to be identical.  But in Theorem~\ref{thm:2}, the
following {\em ad~hoc\/} notion of equivalence will be used.  Note that if
$g(X_0(N^2))=0$, ${\mathfrak h}^{(1)}\sim{\mathfrak h}^{(2)}$ implies
${\mathfrak h}^{(1)}={\mathfrak h}^{(2)}$.

\begin{definition}
If $\mathfrak{h}^{(1)},\mathfrak{h}^{(2)}$ are analytic functions in a
neighborhood of the cusp $\left[\frac1{N^2}\right]$ on~$X_0(N^2)$, i.e., of
the point~$\tau={\rm i}\infty$, $\mathfrak{h}^{(1)}\sim\mathfrak{h}^{(2)}$
signifies that they are the unique (normalized) analytic solutions there of
$L_1\mathfrak{h}=0$, $L_2\mathfrak{h}=0$, two Fuchsian differential
equations satisfying the following conditions.  (i)~Each satisfies the
conditions of Theorem~\ref{thm:ford1}: any ratio of independent solutions
is of the form $(a\tau+b)/(c\tau+d)$, the singular points are the fixed
points of~$\Gamma_0(N^2)$, and the exponent differences are $1/k$ at each
fixed point of order~$k$ and zero at each cusp.  (ii)~The {\em exponents\/}
(not merely exponent differences) of~$L_1,L_2$ are the same.
\end{definition}

\begin{theorem}
\label{thm:2}
Suppose that $X_0(N)$ is of genus zero, and that $f_{N^2}$~is an element of
the function field of~$X_0(N^2)$ that satisfies the divisor condition
\begin{align*}
{\rm div}(f_{N^2})&=
\psi(N)\cdot\left[
({\phi_N'})^*\left(\left[\tfrac11\right]_N/N\right)
-({\phi_N})^*\left(\left[\tfrac11\right]_N\right)
\right],\\
&=\psi(N)\cdot\!\!\!\sum_{\substack{d\ {\rm s.t.}\,d|N|N^2,\\a=a_1,\dots,a_{\varphi(f_{d,N^2})}}}
\Bigl[\left(\left[\tfrac ad\right]_{N^2}\right)-\left(\left[\tfrac11\right]_{N^2}\right)\Bigr],
\end{align*}
where $(\phi_N)^*,(\phi_N')^*$ lift divisors on~$X_0(N),X_0(N)'$ to those
on~$X_0(N^2)$, and where the second equality comes from
Lemma\/~{\rm\ref{lem:cusplifting}}.  Then the hypergeometric `identity'
\begin{equation}
\label{eq:tired}
h_N\left(\phi_N(\cdot)\right) \sim
\left[f_{N^2}(\cdot)/f_{N^2}(\left[\tfrac1{N^2}\right]_{N^2})\right]^{-1/12}
h_N\left(\phi_N'(\cdot)\right)
\end{equation}
holds in a neighborhood of the cusp\/ $\left[\frac1{N^2}\right]$
on~$X_0(N^2)$, i.e., of the point~$\tau={\rm i}\infty$.
\end{theorem}

\begin{proof}
As in the last proof, consider $\tilde{\tilde M}_{N^2},\tilde{\tilde
M}_{N^2}'$, the liftings of~$\tilde{M}_N,\tilde{M}_N'$ to~$X_0(N^2)$
along~$\phi_N,\phi_N'$.  They are projectively the same, in the sense that
any ratio of independent local solutions of $\tilde{\tilde M}_{N^2}\tilde
u=0$, and also of $\tilde{\tilde M}_{N^2}'\tilde u=0$, is of the form
$(a\tau+b)/(c\tau+d)$.  So they have the same exponent differences at each
point of~$X_0(N^2)$.  But they have different exponents.  At~any cusp
of~$X_0(N^2)$ above $x_N=\infty$, resp.~$x_N'=\infty$, i.e., above
$\left[\frac11\right]_N$, resp.~$\left[\frac11\right]_N/N$, the exponents
of $\tilde{\tilde M}_{N^2}$, resp.~$\tilde{\tilde M}_{N^2}'$, are
$n\left(\frac1{12}\psi(N),\frac1{12}\psi(N)\right)$, where $n$~is the cusp
multiplicity.  Converting $\tilde{\tilde M}_{N^2}'$ to a {\em weak\/}
lifting by including a lifting prefactor~$f_{N^2}^{-1/12}$, where
$f_{N^2}$~satisfies the divisor condition of the theorem, will ensure that
$\tilde{\tilde M}_{N^2}$ and~$\tilde{\tilde M}_{N^2}'$ have the same
exponents.  The two sides of~(\ref{eq:tired}) are their respective unique
(normalized) analytic local solutions in a neighborhood of the
cusp~$\left[\frac1{N^2}\right]_{N^2}$.
\end{proof}

The divisor of the theorem is principal, i.e., such a function~$f_{N^2}$
always exists.  The following proposition supplies an explicit formula
for~it.
\begin{proposition}
\label{prop:supplement}
The automorphic function $f_{N^2}:X_0(N^2)\to\mathbb{P}^1({\mathbb{C}})$ of
the previous theorem can be chosen to be
\begin{equation}
\label{eq:cuspprod}
f_{N^2}(\cdot)=
\prod_{\substack{d\ {\rm s.t.}\,1<d|N,\\a=a_1,\dots,a_{\varphi(f_{d,N})}}}
\frac
{\left[x_N(\phi_N(\cdot))-x_N(\left[\frac{a}d\right]_N)\right]^{e_{N/d,N}}}
{\left[x_N'(\phi_N'(\cdot))-x_N(\left[\frac{a}d\right]_N)\right]^{e_{d,N}}},
\end{equation}
a product over the $\sigma_\infty(N)-1$ cusps of~$X_0(N)$ other
than~$\left[\frac11\right]_N$, i.e., than~$\tau=0$.  Here
$e_{d,N}:=N/d(d,N/d)$ is the width of the cusp~$\left[\tfrac{a}d\right]_N$,
as~above.
\end{proposition}
\begin{proof}
On~$X_0(N)$, the function $x_N-x_N\left(\left[\frac{a}d\right]_N\right)$
has divisor
$\left(\left[\frac{a}d\right]_N\right)-\left(\left[\frac11\right]_N\right)$,
the latter term coming from the pole of~$x_N$ at~$\tau=0$.  Similarly
on~$X_0(N)'$, the function $x_N'-x_N\left(\left[\frac{a}d\right]_N\right)$
has divisor
$\left(\left[\frac{a}d\right]_N/N\right)-\left(\left[\frac11\right]_N/N\right)$.
So on~$X_0(N^2)$, $f_{N^2}$~has divisor
\begin{displaymath}
\sum_{\substack{d\ {\rm s.t.}\,d|N,\\a=a_1,\dots,a_{\varphi(f_{d,N})}}}
\!\!\!\!\!
\Bigl[
e_{N/d,N}\cdot (\phi_N)^*\left(\left[\tfrac{a}d\right]_N - \left[\tfrac11\right]_N\right)
-
e_{d,N}\cdot (\phi_N')^*\left(\left[\tfrac{a}d\right]_N/N - \left[\tfrac11\right]_N/N\right)
\Bigr],
\end{displaymath}
where the restriction to $d>1$ has been innocuously dropped.  This
expression splits naturally into two sub-expressions, the first of which is
\begin{subequations}
\begin{align}
\label{eq:fooa}
&\sum_{\substack{d\ {\rm s.t.}\,d|N,\\a=a_1,\dots,a_{\varphi(f_{d,N})}}}
\!\!\!\!\!
\Bigl[
e_{d,N}\cdot (\phi_N')^*\left(\left[\tfrac11\right]_N/N\right)
-
e_{N/d,N}\cdot (\phi_N)^*\left(\left[\tfrac11\right]_N\right)
\Bigr]\\
&\qquad\qquad\qquad\quad
=
\psi(N)\cdot\left[
({\phi_N'})^*\left(\left[\tfrac11\right]_N/N\right)
-({\phi_N})^*\left(\left[\tfrac11\right]_N\right)
\right].
\label{eq:foob}
\end{align}
\end{subequations}
Each cusp of the form~$\left[\frac{a}d\right]_N$ has ramification
index~$e_{d,N}$ over the unique cusp of~$X(1)$, so the sum $\sum e_{d,N}$
over all $\sigma_\infty(N)$~cusps must equal~$\psi(N)$, the degree of the
covering map, as must $\sum e_{N/d,N}$, by symmetry; which explains the
equality between (\ref{eq:fooa}) and~(\ref{eq:foob}).  The
divisor~(\ref{eq:foob}) is identical to the desired divisor of
Theorem~\ref{thm:2}, so it remains to show that the second sub-expression,
namely
\begin{displaymath}
\sum_{\substack{d\ {\rm s.t.}\,d|N,\\a=a_1,\dots,a_{\varphi(f_{d,N})}}}
\!\!\!\!\!
\Bigl[
e_{N/d,N}\cdot (\phi_N)^*\left(\left[\tfrac{a}d\right]_N\right)
\Bigr]
\,\,\,-
\!\!\!\!\sum_{\substack{d\ {\rm s.t.}\,d|N,\\a=a_1,\dots,a_{\varphi(f_{d,N})}}}
\!\!\!\!\!
\Bigl[
e_{d,N}\cdot (\phi_N')^*\left(\left[\tfrac{a}d\right]_N/N\right)
\Bigr],
\end{displaymath}
equals zero.  The first term can be viewed as a sum over the
$\sigma_\infty(N^2)$ cusps of~$X_0(N^2)$.  Any cusp~$\left[\frac{\tilde
a}{\tilde d}\right]_{N^2}$ with $\tilde d|N^2$ that is sent to a cusp of
the form~$\left[\frac{a}d\right]_N$ by~$\phi_N$ will appear in this sum
with multiplicity equal to $e_{N/d,N}=d/(d,N/d)$ times its ramification
index, which by Lemma~\ref{lem:cusplifting2} equals
$e_{d,N^2}/e_{d,N}=N(d,N/d)/d$.  So each of the $\sigma_\infty(N^2)$ cusps
of~$X_0(N^2)$ appears with multiplicity $e_{N/d,N}\,e_{d,N^2}/e_{d,N}=N$.
The subtrahend can be obtained from the first term by applying the Fricke
involution~$w_{N^2}$, which interchanges each cusp of the form
$\left[\frac{\tilde a}{\tilde d}\right]_{N^2}$ with a cusp of the form
$\left[\frac{{\tilde a}'}{N^2/\tilde d}\right]_{N^2}$.  So it equals the
first term, and their difference equals zero.
\end{proof}

\section{Explicit Formulas: The Cases $N=2,\dots,7$}
\label{sec:examples}

\subsection{The $g(X_0(N))=0$, $g(X_0(N^2))=0$ cases ($N=2,3,4,5$)}
\label{subsec:g0}

When $N$ equals $2,3,4$, or~$5$, the fundamental weight-$1$ modular
form~$h_N$ for~$\Gamma_0(N)$ can be expressed in~terms of~${}_2F_1$ on a
neighborhood of the point~$x_N=0$ of~$X_0(N)$ by
\begin{subequations}
\def\theequation{\arabic{section}.\arabic{parentequation}\alph{equation}}
{
\makeatletter
\def\@currentlabel{\arabic{section}.\arabic{parentequation}}
\label{eq:h2def}
\addtocounter{parentequation}{1}
\def\@currentlabel{\arabic{section}.\arabic{parentequation}}
\label{eq:h3def}
\addtocounter{parentequation}{1}
\def\@currentlabel{\arabic{section}.\arabic{parentequation}}
\label{eq:h4def}
\addtocounter{parentequation}{1}
\def\@currentlabel{\arabic{section}.\arabic{parentequation}}
\label{eq:h5def}
\addtocounter{parentequation}{-3}
\makeatother
}
\begin{align}
\label{eq:h2defa}
h_2(z)&=\left[\tfrac1{16^3}(z+16)^3\right]^{-1/12}{}_2F_1\left(\tfrac1{12},\tfrac5{12};\,1;\,\tfrac{1728\,z}{(z+16)^3}\right),\\
\label{eq:h2defb}
&={}_2F_1(\tfrac14,\tfrac14;\,1;\;-z/64);
\\[\jot]
\addtocounter{parentequation}{1}
\setcounter{equation}{0}
\label{eq:h3defa}
h_3(z)&=\left[\tfrac1{{3}^6}(z+3)^3(z+27)\right]^{-1/12}{}_2F_1\left(\tfrac1{12},\tfrac5{12};\,1;\,\tfrac{1728\,z}{(z+3)^3(z+27)}\right),\\
\label{eq:h3defb}
&={}_2F_1(\tfrac13,\tfrac13;\,1;\;-z/27);
\displaybreak[0]\\[\jot]
\addtocounter{parentequation}{1}
\setcounter{equation}{0}
\label{eq:h4defa}
h_4(z)&=\left[\tfrac1{{16}^3}(z^2+16z+16)^3\right]^{-1/12}{}_2F_1\left(\tfrac1{12},\tfrac5{12};\,1;\,\tfrac{1728\,z(z+16)}{(z^2+16z+16)^3}\right),\\
\label{eq:h4defb}
&={}_2F_1(\tfrac12,\tfrac12;\,1;\;-z/16);
\\[\jot]
\addtocounter{parentequation}{1}
\setcounter{equation}{0}
\label{eq:h5defa}
h_5(z)&=\left[\tfrac1{5^3}(z^2+10z+5)^3\right]^{-1/12}{}_2F_1\left(\tfrac1{12},\tfrac5{12};\,1;\,\tfrac{1728\,z}{(z^2+10z+5)^3}\right);
\end{align}
\end{subequations}
where $z$ signifies~$x_N$.  The formulas
(\ref{eq:h2defa}),(\ref{eq:h3defa}),(\ref{eq:h4defa}),(\ref{eq:h5defa})
follow from Definition~\ref{def:hNdef} and the formulas of the form
$j=P_N(x_N)/Q_N(x_N)$ for the covering $\pi_N:X_0(N)\to
X(1)\cong\mathbb{P}^1({\mathbb{C}})\ni j$ given in the appendix.  The
simpler expressions (\ref{eq:h2defb}),(\ref{eq:h3defb}),(\ref{eq:h4defb})
follow from the fact mentioned after Proposition~\ref{prop:2}: each of
$h_2,h_3,h_4$ satisfies an equation of hypergeometric type.  For instance,
lifting ${}_2\mathcal{E}_1(\frac12,\frac5{12};1)$ from the $\hat
J$-sphere~$X(1)$ to~$X_0(2)\ni x_2=:z$ yields
\begin{equation}
\left\{D_z^2 + \left[\frac1z+\frac{1}{2(z+64)}\right]D_z + 
\left[
\frac{1}{16\,z(z+64)}
\right]\right\}h_2=0,
\end{equation}
This is reduced to the Gau\ss\ equation
${}_2\mathcal{E}_1(\frac14,\frac14;1)$ by the linear scaling $\hat
z=-z/64$, yielding~(\ref{eq:h2defb}).  The expressions
(\ref{eq:h3defb}),(\ref{eq:h4defb}) follow likewise.  The equivalence of
the `a'~and~`b' expressions yields a cubic, a quartic, and a sextic
transformation of~${}_2F_1$.  These turn~out to be special cases of three
of Goursat's hypergeometric transformations, listed as
(\oldstylenums{121}),(\oldstylenums{129}),(\oldstylenums{135}) in his
classical tabulation~\cite{Goursat1881}.

There is no expression~(\ref{eq:h5def}b), since the Fuchsian equation
obtained by lifting ${}_2\mathcal{E}_1(\frac12,\frac5{12};1)$ to~$X_0(5)\ni
x_5=:z$, which appears as~(\ref{eq:h5de}) above, has four singular points
rather than three.  Its analytic solution~$h_5$ near~$z=0$ can be reduced
to~$\Hl$ by an appropriate M\"obius transformation $\hat z=(Az+B)/(Cz+D)$.
By definition, $\Hl(a,q;\alpha,\beta,\gamma,\delta;\cdot)$ is the unique
(normalized) analytic local solution at~$\hat z=0$ of a canonical Fuchsian
equation with the four singular points $\hat z=0,1,a,\infty$.  (The
parameters $\alpha,\beta,\gamma,\delta$ are exponent-related; $q$~is the
accessory parameter.)  So the transformation
$\left(\begin{smallmatrix}A&B\\C&D\end{smallmatrix}\right)\in{\it
PGL}(2,\mathbb{C})$ must take the singular points $0,-11\pm2{\rm i},\infty$
of~(\ref{eq:h5de}) to~$0,1,a,\infty$ for
some~$a\in\mathbb{C}\setminus\{0,1\}$.  (The points $x_5=-11\pm2{\rm i}$
are not cusps.  They are the order-$2$ elliptic fixed points of~$X_0(5)$,
arising from self-isogenies of lemniscatic elliptic curves;
cf.~\cite[p.~48]{Elkies98}.)  The resulting alternative expression
for~$h_5$ is not pleasing, since it involves the radical ${\rm
i}=\sqrt{-1}$, and is not given here.

\begin{proposition}
Let\/ $\mathbb{C}$-valued analytic functions $h_2,h_3,h_4$ be defined in a
neighborhood of\/ $0\in\mathbb{C}$
by\/~{\rm(\ref{eq:h2def}),(\ref{eq:h3def}),(\ref{eq:h4def})}, in~terms
of~${}_2F_1$.  Then for all~$x$ near\/~$0$,
\begin{align*}
&h_2\left(x(x+16)\right)\\[-\jot]
&\qquad=2\left[x+16\right]^{-1/4} 
  h_2\left(\frac{x^2}{x+16}\right),\displaybreak[0]\\
&h_3\left(x(x^2+9x+27)\right)\\[-\jot]
&\qquad=3\left[x^2+9x+27\right]^{-1/3}
  h_3\left(\frac{x^3}{x^2+9x+27}\right),\\
&h_4\left(x(x+4)(x^2+4x+8)\right)\\[-\jot]
&\qquad=4\left[(x+2)(x^2+4x+8)\right]^{-1/2}
  h_4\left(\frac{x^4}{(x+2)(x^2+4x+8)}\right).
\end{align*}
That is, they respectively satisfy quadratic, cubic, and quartic functional
equations.
\end{proposition}
\begin{proof}
These follow from Theorem~\ref{thm:1} and the appendix.  The argument of
the left~$h_N$ is the polynomial expression $R(x_{N^2})$ for~$x_N$ in~terms
of~$x_{N^2}$, and that of the right~$h_N$ is the rational expression
$x_{N^2}^N/S'(x_{N^2})$ for~$x_N'$ in~terms of~$x_{N^2}$.
\end{proof}

The corresponding identity involving~$h_5$ appeared above as
Proposition~\ref{prop:h5} and will not be repeated here.  The three-term
recurrence~(\ref{eq:h5recurrence}), which is a simpler characterization
of~$h_5$ than any expression involving ${}_2F_1$ or~$\Hl$, is obtained by
substituting $h_5(z)=\sum_{n=0}^\infty c_nz^n$ into the lifted differential
equation~(\ref{eq:h5de}).

\begin{corollary}
The following hypergeometric identities are valid 
for all~$\hat x$ in a neighborhood of\/~$0$.
\begin{align*}
{}_2F_1\left(\tfrac14,\tfrac34;\,1;\,1-\Bigl(\frac{1-\hat x}{1+3\hat x}\Bigr)^2\right)
&=
(1+3\hat x)^{1/2}\, {}_2F_1\left(\tfrac14,\tfrac34;\,1;\,\hat x^2\right),\\
{}_2F_1\left(\tfrac13,\tfrac23;\,1;\,1-\Bigl(\frac{1-\hat x}{1+2\hat x}\Bigr)^3\right)
&=
(1+2\hat x)\,\, {}_2F_1\left(\tfrac13,\tfrac23;\,1;\,\hat x^3\right),\\
{}_2F_1\left(\tfrac12,\tfrac12;\,1;\,1-\Bigl(\frac{1-\hat x}{1+\hat x}\Bigr)^4\right)
&=
(1+\hat x)^2\, {}_2F_1\left(\tfrac12,\tfrac12;\,1;\,\hat x^4\right).
\end{align*}
\end{corollary}
\begin{proof}
(i)~Express $h_2,h_3,h_4$ above in~terms of~${}_2F_1$.  (ii)~Apply Pfaff's
transformation
${}_2F_1(\alpha,\beta;\gamma;s)=(1-s)^{-\alpha}{}_2F_1(\alpha,\gamma-\beta;\gamma;s/(s-1))$
to each~${}_2F_1$.  (iii)~Rewrite each equation in~terms of a new
(M\"obius-transformed) Hauptmodul~$\hat x$, defined to equal
$x/(x+32),x/(x+9),x/(x+4)$, respectively.
\end{proof}

These three transformations of~${}_2F_1$ are the quadratic
arithmetic--geometric mean (AGM) iteration in signature~$4$, the cubic one
in signature~$3$, and the quartic one in signature~$2$.  The most familiar
is the last, which is classical.  It~is an iterate of Landen's quadratic
transformation of the first complete elliptic integral function
$\mathsf{K}(\cdot)=\frac\pi2\,{}_2F_1(\frac12,\frac12;1;\cdot)$.  The other
two were found by Ramanujan~\cite[pp.~97, 146]{BerndtV}.  They play a major
role in his theory of elliptic functions to alternative bases.

Our new derivation of these hypergeometric transformations is of~interest
not only because it clarifies their modular origin, but because it extends
to~$N=5,6,7$.  Proposition~\ref{prop:h5}, the $N=5$ transformation, does
not fit into Ramanujan's framework.  The $N=6,7$ transformations are even
more exotic.

\subsection{The $g(X_0(N))=0$, $g(X_0(N^2))>0$, $g(X_0^+(N^2))=0$ cases ($N=6,7$)}
\label{subsec:g1}

When $N$ equals $6$ or~$7$, the fundamental weight-$1$ modular form~$h_N$
for~$\Gamma_0(N)$ can be expressed in~terms of~${}_2F_1$ on a neighborhood
of the point~$x_N=0$ of~$X_0(N)$ by
\begin{subequations}
\def\theequation{\arabic{section}.\arabic{parentequation}\alph{equation}}
{
\makeatletter
\def\@currentlabel{\arabic{section}.\arabic{parentequation}}
\label{eq:h6def}
\addtocounter{parentequation}{1}
\def\@currentlabel{\arabic{section}.\arabic{parentequation}}
\label{eq:h7def}
\addtocounter{parentequation}{-1}
\makeatother
}
\begin{align}
\label{eq:h6defa}
\begin{split}
h_6(z)&=\left[\tfrac{1}{2^{12}3^6}(z+6)^3(z^3+18z^2+84z+24)^3\right]^{-1/12}\\
&\qquad\qquad\qquad\times{}_2F_1\left(\tfrac1{12},\tfrac5{12};\,1;\,
\tfrac{1728\,z(z+9)^2(z+8)^3}{(z+6)^3(z^3+18z^2+84z+24)^3}\right)
\end{split}
\\
\label{eq:h6defb}
&=\Hl\,(\tfrac98,\tfrac34;\,1,1,1,1;\,-z/8),
\displaybreak[0]\\[\jot]
\addtocounter{parentequation}{1}
\setcounter{equation}{0}
\label{eq:h7defa}
\begin{split}
h_7(z)&=\left[\tfrac{1}{49}(z^2+5z+1)^3(z^2+13z+49)\right]^{-1/12}\\ 
&\qquad\times{}_2F_1\left(\tfrac1{12},\tfrac5{12};\,1;\,
\tfrac{1728\,z}{(z^2+5z+1)^3(z^2+13z+49)}\right);
\end{split}
\end{align}
\end{subequations}
where $z$ signifies~$x_N$.  The formulas
(\ref{eq:h6defa}),(\ref{eq:h7defa}) follow from Definition~\ref{def:hNdef},
together with the formulas $j=P_N(x_N)/Q_N(x_N)$ for the covering
$\pi_N:X_0(N)\to X(1)$ given in the appendix.  Lifting
${}_2\mathcal{E}_1(\frac12,\frac5{12};1)$ from the $\hat J$-sphere~$X(1)$
to~$X_0(N)\ni x_N=:z$ yields an equation based on an operator
like~(\ref{eq:Poole}), which if $N=6,7$ is respectively
\begin{align}
\label{eq:h6de}
&\left\{D_z^2 + \left[\frac1z+\frac1{z+8}+\frac1{z+9}\right]D_z + 
\left[
\frac{z+6}{z(z+8)(z+9)}
\right]\right\}h_6=0,\\
\label{eq:h7de}
&\left\{D_z^2 + \left[\frac1z+\frac{2(2\,z+13)}{3(z^2+13\,z+49)}\right]D_z + 
\left[
\frac{4\,z+21}{9\,z(z^2+13\,z+49)}
\right]\right\}h_7=0.
\end{align}
The four singular points of~(\ref{eq:h6de}) are $0,-8,-9,\infty$, so a
linear rescaling $\hat z=-z/8$ yields the much simpler expression for~$h_6$
given in~(\ref{eq:h6defb}), in~terms of~$\Hl$.  The singular points
of~(\ref{eq:h7de}) are $0,-\frac{13}2\pm\frac32\sqrt{-3},\infty$.  (The
points $x_7=-\frac{13}2\pm\frac32\sqrt{-3}$ are not cusps.  They are the
order-$3$ elliptic fixed points of~$X_0(7)$, arising from self-isogenies of
equianharmonic elliptic curves.)  So when $N=7$, just as when $N=5$, the
three finite singular points of the lifted equation are not collinear.
A~reduction of~$h_7$ to~$\Hl$ must accordingly be performed by an awkward
M\"obius transformation that involves the radical~$\sqrt{-3}$.  The
resulting expression is unpleasing and is omitted.

The $N=6$ case is a bit special, since each of the four singular points
of~(\ref{eq:h6de}) on~$X_0(6)$ is a cusp, there being no~elliptic fixed
points.  That~is, $\sigma_\infty(6)=4$ and $\varepsilon_{\rm
i}(6)=\varepsilon_{\rho}(6)=0$.  It is a result of
Beauville~\cite{Beauville82} and Sebbar~\cite{Sebbar2002} that there are
only six genus-$0$ algebraic curves of the form
$\Gamma_1\setminus\mathfrak{H}^*$, where $\Gamma_1<\Gamma$ is a congruence
group, with exactly four cusps and no~elliptic points.  The number of
curves shrinks to four if isomorphic curves arising from
subgroups~$\Gamma_1$ that are conjugate in~${\it PSL}(2,\mathbb{R})$ are
identified~\cite{Beukers2002,Sebbar2002}.  The list includes~$X_0(6)$, and
also $X_0(8)$, $X_0(9)$, and~$X_1(5)$.

\begin{proposition}
\label{prop:6}
Let $h_6$, a $\mathbb{C}$-valued function, be defined in a neighborhood of
$0\in\mathbb{C}$ by $h_6(z)=\sum_{n=0}^\infty c_nz^n$, where the
coefficients satisfy the three-term recurrence
\begin{equation}
\label{eq:h6recurrence}
n^2\,c_{n-1} + (17n^2+17n+6)\,c_n + 72(n+1)^2\,c_{n+1}=0,
\end{equation}
initialized by $c_{-1}=0$ and $c_0=1$.  Equivalently, let $h_6$ be defined
in~terms of ${}_2F_1$ or~$\Hl$ by\/~{\rm(\ref{eq:h6def}).}  Let algebraic
functions $z=z(t)$, $z'=z'(t)$ be defined implicitly, in a neighborhood of
$t=\infty\in\mathbb{P}^1(\mathbb{C})$, by
\begin{subequations}
\begin{align}
z+72/z' &= (t-2)(t^5-10\,t^4+28\,t^3-26\,t^2+20\,t+4),\\
z\cdot[72/z'\,] &= 72\,t\,(t-2)^2(t^2-t+1).
\end{align}
\end{subequations}
{\rm(}Of~the two branches, the one on which $z,z'\to0$ as~$t\to\infty$ is
to be chosen.{\rm)} Then for all~$t$ in a neighborhood of~$\infty$,
\begin{equation}
\label{eq:sextic}
h_6(z(t)) = 6\left\{\frac{z(t)^6}{z'(t)}
\cdot\frac{[z(t)+9]^3}{[z'(t)+9]^2}
\cdot\frac{[z(t)+8]^2}{[z'(t)+8]^3}
\right
\}^{-1/12} \!h_6(z'(t)).
\end{equation}
\end{proposition}
\begin{remarkaftertheorem}
By examination, $z\sim 72/t$ and $z'\sim 72/t^6$ as~$t\to\infty$.  So
(\ref{eq:sextic})~is a sextic functional equation for~$h_6$.  The
multiplier, i.e., the prefactor multiplying~$h_6$ on the right-hand side,
tends to unity as~$t\to\infty$.
\end{remarkaftertheorem}
\begin{remarkaftertheorem}
The sequence $d_n:=72^nc_n$, $n\ge0$, of Maclaurin coefficients of
$h_6(72z)$ is an integral sequence.  It~begins $1$, $-6$, $42$, $-312$,
$2394$, $-18756$, $149136$, $-1199232$, $\dots$.  A~three-term recurrence
for~$d_n$ equivalent to~(\ref{eq:h6recurrence}) was discovered by
Coster~\cite{Coster83}, in an investigation of Beauville's six curves.
The~representation $d_n=\sum_{k=0}^n
\left(\begin{smallmatrix}n\\k\end{smallmatrix}\right)(-8)^k\sum_{j=0}^{n-k}\left(\begin{smallmatrix}n-k\\j\end{smallmatrix}\right)^3$
was later worked~out by Verrill~\cite[Table~2]{Verrill99}, who began with
the differential equation satisfied by $h_6=[1]^6[6]\,/\,[2]^3[3]^2$ as a
function of $x_6/72=[2][6]^5/\,[1]^5[3]$, i.e., in~effect
with~(\ref{eq:h6de}).  This sequence now appears in Sloane's On-Line
Encyclopedia~\cite{Sloane2005} as sequence~{\tt A093388}.

Perhaps due~to the simplicity of the Heun representation~(\ref{eq:h6defb}),
the function~$h_6$ is combinatorially significant.  For example, the
perimeter generating function for three-dimensional staircase polygons can
be expressed in~terms of~$h_6$~\cite{Guttmann93}.
\end{remarkaftertheorem}
\begin{proof}[Proof of Proposition\/~{\rm\ref{prop:6}}]
The recurrence~(\ref{eq:h6recurrence}) comes by substituting
$h_6(z)=\sum_{n=0}^\infty c_nz^n$ into~(\ref{eq:h6de}).  Otherwise, this
follows from Theorem~\ref{thm:2} and the prefactor formula given in
Proposition~\ref{prop:supplement}, with the understanding that $t,z,z'$
signify $t_{36},x_6,x_6'$, the Hauptmoduln of~$X_0^+(36),X_0(6),X_0(6)'$.
As~the appendix explains, the relation between the sum $x_6+72/x_6'$ and
product $x_6\cdot[72/x_6']$ is parametrized by~$t_{36}$.  The quantity in
curly braces in~(\ref{eq:sextic}) is a product over the
$\sigma_\infty(6)-1=3$ cusps of~$X_0(6)$ other than~$\left[\frac11\right]$.
In~all, the cusps are
$\left[\frac11\right]\ni0,\left[\frac12\right],\left[\frac13\right],\left[\frac16\right]\ni{\rm
i}\infty$, with respective widths (i.e., multiplicities over~$X(1)$) equal
to $e_{1,6}=6$, $e_{2,6}=3$, $e_{3,6}=2$, and~$e_{6,6}=1$.  A~comparison
with the formula~(\ref{eq:12sheeted}) for the covering $j=j(x_6)$ reveals
that these cusps are located respectively at the points
$x_6=\infty,-8,-9,0$ on~$X_0(6)$.  The prefactor follows.

However, Theorem~\ref{thm:2} does not quite imply that the functions
of~$t:=t_{36}$ on the two sides of~(\ref{eq:sextic}) are equal.  What it
implies is that they are the unique (normalized) analytic solutions near
the point $(t_{36},s_{36})=(\infty,-\infty)$ on the elliptic curve
$X_0(36)$ modeled by~(\ref{eq:X036eqn}), i.e., by
\begin{equation}
s_{36}^2 = t_{36}^4 - 8\,t_{36}^3 + 12\,t_{36}^2 - 8\,t_{36} + 4,
\end{equation}
of a pair of differential equations that are the same {\em up~to one degree
of freedom\/}: their single affine accessory parameter.  To
verify~(\ref{eq:sextic}), one must show the two equations are identical.
They are obtained from~(\ref{eq:h6de}), the equation on~$X_0(6)\ni x_6$
satisfied by~$h_6$, by (i)~lifting along $(t_{36},s_{36})\mapsto x_6$, and
(ii)~lifting along $(t_{36},s_{36})\mapsto x'_6$ and including the
prefactor.  Each lifted equation on~$X_0(36)\ni(t_{36},s_{36})=:(t,s)$ is
based on a Fuchsian operator of the normal form~(\ref{eq:Pooleplus}), with
$g=1$ and $n=\sigma_\infty(36)=12$ singular points (i.e., cusps).
Fortunately, the two turn~out to be the same.  Regardless of which route is
taken, the lifting of~(\ref{eq:h6de}) to~$X_0(36)$ is
\begin{align}
\label{eq:foo36}
&\Biggl\{
\left(s\frac{d}{dt}\right)^2
+\left[\frac{s}{t} + \frac{s}{t-1} + \frac{s}{t-2} + \frac{s\,(2t-1)}{t^2-t+1}
+6t+\fbox{$-6$}\,\right]
\left(s\frac{d}{dt}\right)
\\
&\quad\qquad{}+3\left[
\frac{-(s+2)}{t} + \frac{s-6}{t-2} + \frac{-(t+1)s + 3t}{t^2-t+1} + (6t^2 -
23t+8) + 6s\,
\right]\Biggr\}\,\mathfrak{h}_{36}=0.
\nonumber
\end{align}
Here the single affine accessory parameter is boxed, and the $n-3+3g=12$
projective accessory parameters are contained in the second bracketed
expression.

This equation on the equianharmonic elliptic curve $X_0(36)$, every ratio
of solutions of~which is of the form $(a\tau+b)/(c\tau+d)$, is of
independent interest.  It is a uniformizing differential equation of the
sort guaranteed to exist by Theorem~\ref{thm:ford1}, but few such equations
in the case of positive genus have appeared in the literature.  By
examination, it displays the $12$~cusps of~$X_0(36)$ explicitly.  They come
in pairs, each pair being related by the Fricke involution $s_{36}\mapsto
-s_{36}$, and are located at
$t_{36}=0,1,2,\frac12\pm\frac{\sqrt{-3}}2,\infty$.  The $6$~{\em
rational\/} cusps $\left[\frac{a}d\right]_{36}$ are those with
$t_{36}=0,1,\infty$.  They are singled~out by
$\varphi\left((d,N/d)\right)\le2$, a~standard arithmetic rationality
condition~\cite{Ogg73}, and are respectively
$\left[\frac12\right],\left[\frac1{18}\right];\left[\frac1{4}\right],\left[\frac1{9}\right];\left[\frac1{1}\right]\ni0,\left[\frac1{36}\right]\ni{\rm
i}\infty$.

For present purposes, all that matters is that the value of the single
(boxed) affine accessory parameter is independent of which of the two
liftings is used.  It~is actually possible to verify this without deriving
(\ref{eq:foo36}) in~full, by simply working~out the first two terms in an
expansion of its coefficient of $s\,d/dt$ about~$(t,s)=(\infty,-\infty)$,
i.e., about $\tau={\rm i}\infty$.  This can readily be done by~hand, though
why the result is the same for the lifting from $X_0(6)$ to~$X_0(36)$ and
the weak lifting from $X_0(6)'$ to~$X_0(36)$ is not entirely clear.
\end{proof}

\begin{proposition}
\label{prop:7}
Let $h_7$, a $\mathbb{C}$-valued function, be defined in a neighborhood of
$0\in\mathbb{C}$ by $h_7(z)=\sum_{n=0}^\infty c_nz^n$, where the
coefficients satisfy the three-term recurrence
\begin{equation}
\label{eq:h7recurrence}
(3n-1)^2\,c_{n-1} + 3(39n^2+26n+7)\,c_n + 441(n+1)^2\,c_{n+1}=0,
\end{equation}
initialized by $c_{-1}=0$ and $c_0=1$.  Equivalently, let $h_7$ be defined
in~terms of ${}_2F_1$ by~{\rm(\ref{eq:h7defa}).}  Let algebraic functions
$z=z(t)$, $z'=z'(t)$ be defined implicitly, in a neighborhood of
$t=\infty\in\mathbb{P}^1(\mathbb{C})$, by
\begin{subequations}
\begin{align}
z+49/z' &= (t^3-7\,t^2+14\,t-7)(t^4-14\,t^3+63\,t^2-98\,t+35),\\
z\cdot[49/z'\,] &= 49\,(t^3-7\,t^2+14\,t-7)^2.
\end{align}
\end{subequations}
{\rm(}Of~the two branches, the one on which $z,z'\to0$ as~$t\to\infty$ is
to be chosen.{\rm)} Then for all~$t$ in a neighborhood of~$\infty$,
\begin{equation}
\label{eq:septic}
h_7\left(z(t)\right) = 7\left[\frac{z(t)^7}{z'(t)}\right]^{-1/12} \!h_7\left(z'(t)\right).
\end{equation}
\end{proposition}
\begin{remarkaftertheorem}
By examination, $z\sim 49/t$ and $z'\sim 49/t^7$ as~$t\to\infty$.  So
(\ref{eq:septic})~is a septic functional equation for~$h_7$.  The
multiplier tends to unity as~$t\to\infty$.
\end{remarkaftertheorem}
\begin{remarkaftertheorem}
The sequence $d_n:=441^nc_n$, $n\ge0$, of Maclaurin coefficients of
$h_7(441z)$ is an integral sequence.  It~begins $1$, $-21$, $693$, $-23940$,
$734643$, $-13697019$, $-494620749$, $83079255420$, $\dots$.
\end{remarkaftertheorem}
\begin{proof}[Proof of Proposition\/~{\rm\ref{prop:7}}]
This resembles the proof of Proposition~\ref{prop:6}, but is simpler
because $7$~is not composite.  Here $t,z,z'$ are $t_{49},x_7,x_7'$, the
Hauptmoduln of $X_0^+(49),X_0(7),X_0(7)'$.  There are only two cusps
on~$X_0(7)$, namely $\left[\frac11\right]\ni0$,
$\left[\frac17\right]\ni{\rm i}\infty$, with respective widths $e_{1,7}=7$
and~$e_{7,7}=1$.  They are located at~$x_7=\infty,0$, respectively.  So the
product~(\ref{eq:cuspprod}) comprises only a single factor,
namely~$x_7^7/x_7'$, i.e.,~$z(t)^7/z'(t)$.  As in the $N=6$ case, the
lifting of the differential equation~(\ref{eq:h7de}) from $X_0(7)$ to the
elliptic curve $X_0(49)$, and its weak lifting from $X_0(7)'$ to~$X_0(49)$,
turn~out to be identical due~to the values of their affine accessory
parameters being the same.  So the two sides of~(\ref{eq:septic}) are
equal, as claimed.
\end{proof}

The cases $N=8,9$ can in principle also be handled by Theorem~\ref{thm:2},
like~${N=6,7}$, since equations for the curves $X_0(64),X_0(81)$ are
known~\cite{Shimura95}.  ($X_0(64),X_0(81)$~are non-hyperelliptic of
genera~$3,4$, and $X_0(64)$ is the degree-$4$ Fermat curve.)  The covering
maps $\phi_8,\phi_8',\phi_9,\phi_9'$ could also be worked~out.  But for
both $N=8$ and~$9$, a simple parametrization of the relation between $x_N$
and~$x_N'$ may be lacking, since $X_0^+(64),X_0^+(81)$ are of genus~$1$,
not~$0$.

\section{Possible Extensions}
\label{sec:extensions}

We conclude by mentioning some possible generalizations of our approach.
Each identity derived above arises from the covering $\phi_N:X_0(N^2)\to
X_0(N)$ induced by a subgroup relation $\Gamma_0(N)>\Gamma_0(N^2)$.  Under
the Fricke involution~$w_{N^2}$, $\Gamma_0(N)$~is conjugated to
$\Gamma_0(N)':=w_{N^2}^{-1} \Gamma_0(N) w_{N^2}>\Gamma_0(N^2)$, and
$\phi_N':X_0(N^2)\to X_0(N)'$ arises from the latter relation.  As a first
generalization, one could begin instead with the covering
$\phi_{N,M}:X_0(MN)\to X_0(N)$ induced by $\Gamma_0(N)>\Gamma_0(MN)$, where
$M\neq N$ is allowed.  The canonical Hauptmoduln of $X_0(MN),X_0(N)$ are
$x_{MN},x_N$, and that of $w_{MN}^{-1}\Gamma_0(N)w_{MN}>\Gamma_0(MN)$ is
$x_N''=x_N''(\tau):=x_N(M\tau)$.  If for~example $N=2$ and~$M=3$, this
approach turns~out to yield the functional equation
\begin{equation}
h_2\left(\frac{x(x+8)^3}{x+9}\right)
= 3\left[x+9\right]^{-1/2} h_2\left(\frac{x^3(x+8)}{(x+9)^3}\right),
\end{equation}
where $x:=x_6$.  The argument of the left~$h_2$ is the
expression~(\ref{eq:x2x3}a) for $x_2$ in~terms of~$x_6$, and that of the
right is $x_2''=x_2(3\cdot)$ in~terms of~$x_6$.  This can be rewritten as
\begin{multline}
{}_2F_1\left(\tfrac14,\tfrac34;\,1;\, \frac{64\,p}{(3+6p-p^2)^2}  \right)\\
=\left(\frac{9(3+6p-p^2)}{27-18p-p^2}\right)^{1/2} \!{}_2F_1\left(\tfrac14,\tfrac34;\,1;\,\frac{64\,p^3}{(27-18p-p^2)^2}\right),
\end{multline}
a hypergeometric identity valid near $p=0$.  This is a known identity: the
cubic modular equation in Ramanujan's theory of elliptic functions in
signature~$4$, which has been proved by other
means~\cite[pp.~152--3]{BerndtV}.  The parameter~$p$ is the
M\"obius-transformed Hauptmodul $x_6/(x_6+8)=9\,[1]^4[6]^8/\,[2]^8[3]^4$
of~$X_0(6)$, which has a zero at the cusp $\left[\tfrac16\right]_6\ni{\rm
i}\infty$ and a pole at the cusp~$\left[\tfrac12\right]_6$.  A~second
noteworthy functional equation comes from the choices $N=5$ and~$M=2$.
It~is
\begin{equation}
\label{eq:lastident}
h_5\left(\frac{x(x+5)^2}{x+4}\right)
= 2\left[x+4\right]^{-1/2} h_5\left(\frac{x^2(x+5)}{(x+4)^2}\right),
\end{equation}
where $x:=x_{10}$.  The argument of the left~$h_5$ is the expression for
$x_5$ in~terms of~$x_{10}$, and that of the right is $x_5''=x_5(2\cdot)$
in~terms of~$x_{10}$, both expressions being due~to Fricke.  Here
$x_{10}$~is the canonical Hauptmodul $20\,[2][10]^3/\,[1]^3[5]$
of~$X_0(10)$.  The identity~(\ref{eq:lastident}) merits comparison with
Proposition~\ref{prop:h5}.

A further generalization comes from focusing on~$h_N$, originally given in
terms of~${}_2F_1$, as a weight-$1$ modular form for~$\Gamma_0(N)$ that
necessarily satisfies a second-order differential equation on~$X_0(N)$.
As~an abstraction of this situation, let $f$ be any weight-$1$ modular form
for a congruence subgroup~$G$ of genus zero, with
Hauptmodul~$x_G:=x_G(\tau)$.  Any $g\in{\it SL}(2,\mathbb{R})$ will yield a
projective action on~$\mathfrak{H}\ni\tau$ and a weight-$1$ modular form
$f':=f(x_G(g(\tau)))$ for the conjugated group $G':=g^{-1}Gg$, which has
Hauptmodul $x_{G'}=x_{G'}(\tau):=x_G(g(\tau))$.  So both $f$ and~$f'$ will
be weight-$1$ modular forms for $\tilde G:=G\cap G'$.  The simplest case is
when $\tilde G$~is of genus zero as~well, with a Hauptmodul~$x_{\tilde G}$.
In~this case $x_G$,~$x_{G'}$, and the weight-zero form~$f'/f$ will all be
expressible in~terms of~$x_{\tilde G}$, yielding an identity of the form
\begin{equation}
\label{eq:last}
f(x_G(x_{\tilde G})) = A(x_{\tilde G})\,f(x_{G'}(x_{\tilde G})),
\end{equation}
with $A(\cdot)$ a fractional power of some rational function.  The functional
equations derived above for~$h_N$, $N=2,\dots,5$, were all of this type.

In general the identity~(\ref{eq:last}) will have no interpretation as a
transformation of~${}_2F_1$, since the differential equation satisfied
by~$f$ will not be of hypergeometric type (though it will be of second
order).  However, an additional extension suggests itself.  If $f$~has
weight~$k$, $k\ge1$, the differential equation will be of order~$k+1$.  For
example, it could be the order-$(k+1)$ differential equation satisfied by
any generalized hypergeometric function of the type~${}_{k+1}F_k$.  The
possibility therefore exists that by an approach similar to the one used in
this article, one could derive functional equations satisfied by
certain~${}_{k+1}F_k$, $k>1$.  This remains to be explored.

\appendix
\section*{Appendix. Hauptmoduln}
\renewcommand\thesection{A}

If $X_0(N)$ is of genus~$0$ then it has a Hauptmodul~$x_N$, which may be
regarded as a function on~$\mathfrak{H}\ni\tau$.  This Hauptmodul is
determined uniquely by the condition that it have divisor $({\rm
i}\infty)-(0)$, i.e., a simple zero (resp.~pole) at the cusp $\tau={\rm
i}\infty$ (resp.~the cusp~$\tau=0$), together with a standard normalization
condition that $x_N(-N/\tau)$ have a Fourier expansion on~$\mathfrak{H}$
that begins $q^{-1}+\cdots$, where $q:=e^{2\pi{\rm i}\tau}$.  The
$j$-invariant will be a rational function of~$x_N$, and $x_N$~itself will
be in the function field of~$X_0(N^2)$, and in~fact of~$X_0(MN)$ for
every~$M\ge1$.  

This appendix collects the formulas of the form $j=j(x_N)$, and for $x_N$
itself, which are used in the body of the article.  N.~Fine's notation
$[k]$ for the function $\eta(k\cdot)$, i.e., for $\eta(k\tau)$,
$\tau\in\mathfrak{H}$, is used to save space.  Here $\eta$~is the Dedekind
eta function $q^{1/24}\prod_{n\in\mathbb{N}}(1-q^n)$, satisfying
$\eta(\tau+1)=e^{\pi{\rm i}/12}\eta(\tau)$, $\eta(-1/\tau)=(-{\rm
i}\tau)^{1/2}\eta(\tau)$.

The simplest cases are $N=2,3,4,5$, when $X_0(N^2)$ like~$X_0(N)$ is of
genus~$0$, and a Hauptmodul~$x_{N^2}$ exists.

\begin{itemize}
\item $N=2$.  $x_2:=2^{12}\,[2]^{24}/\,[1]^{24}$, and $x_2$~satisfies
$x_2(\tau)\,x_2(-1/2\tau)=2^{12}$.  The equality $j=(x_2+16)^3/x_2$ holds.
$x_4:=2^8\,[4]^8/\,[1]^8$ is the Hauptmodul on~$X_0(4)$, and
$x_2=x_4(x_4+16)$.  The Hauptmodul~$x_4$ satisfies
$x_4(\tau)\,x_4(-1/4\tau)=2^8$.
\item $N=3$.  $x_3:=3^6\,[3]^{12}/\,[1]^{12}$, and $x_3$~satisfies
$x_3(\tau)\,x_3(-1/3\tau)=3^6$.  The equality $j=(x_3+3)^3(x_3+27)/x_3$
holds.  $x_9:=3^3\,[9]^3/\,[1]^3$ is the Hauptmodul on~$X_0(9)$, and
$x_3=x_9(x_9^2+9x_9+27)$.  The Hauptmodul~$x_9$ satisfies
$x_9(\tau)\,x_9(-1/9\tau)=3^3$.
\item $N=4$. $x_4:=2^8\,[4]^8/\,[1]^8$, and $x_4$~satisfies
$x_4(\tau)\,x_4(-1/4\tau)=2^8$.  The equality
$j=(x_4^2+16x_4+16)^3/\left(x_4(x_4+16)\right)$ holds.
$x_{16}:=2^3\,[2][16]^2/\,[1]^2[8]$ is the Hauptmodul on~$X_0(16)$, and
$x_4=x_{16}(x_{16}+4)(x_{16}^2+4x_{16}+8)$.  The Hauptmodul~$x_{16}$
satisfies $x_{16}(\tau)\,x_{16}(-1/16\tau)=2^3$.
\item $N=5$. $x_5:=5^3\,[5]^6/\,[1]^6$, and $x_5$~satisfies
$x_5(\tau)\,x_5(-1/5\tau)=5^3$.  The equality $j=(x_5^2+10x_5+5)^3/x_5$
holds.  $x_{25}:=5\,\,[25]\,/\,[1]$ is the Hauptmodul on~$X_0(25)$, and
$x_5=x_{25}(x_{25}^4+5x_{25}^3+15x_{25}^2+25x_{25}+25)$.  The
Hauptmodul~$x_{25}$ satisfies $x_{25}(\tau)\,x_{25}(-1/25\tau)=5$.
\end{itemize}

\noindent
Of the preceding formulas, the most familiar are the definition of~$x_N$ as
an eta product and the expression for~$j$ in~terms of~$x_N$, when
$N=2,3,5$.  Like the corresponding formulas for $N=7,13$, they go back to
Weber, indeed to Klein, and are often reproduced~\cite[\S\,4]{Elkies98}.
The remaining formulas were extracted from Fricke~\cite[II.~Abschnitt,
4.~Kap.]{Fricke22}, with due care.  Fricke occasionally varied the
normalization of his Hauptmoduln~$\tau_N$, corresponding to~$x_N$, and his
$\tau_2(\cdot)$ is proportional to~$x_2(2\cdot)$ rather than
to~$x_2(\cdot)$.  He also did not give eta-product expressions for his
$\tau_4,\tau_{16}$, but they are easily worked~out~\cite{Newman58}.
Knopp~\cite[\S\,7.6]{Knopp70} gives an independent derivation of the
quintic formula for $x_5$ in~terms of~$x_{25}$.

From each expression for $j$ in~terms of~$x_N$, one obtains an expression
for $j_N$, where $j_N(\tau):=j(N\tau)$, by applying the involution
$\tau\mapsto-1/\tau$ to both sides.  Thus, $j_2=(x_2+256)^3/x_2^2$.
Similarly, applying the Fricke involution $w_{N^2}:\tau\mapsto-1/N^2\tau$
to both sides of the polynomial expression for $x_N$ in~terms of~$x_{N^2}$
yields a non-polynomial rational expression for~$x_N'$, where
$x_N'(\tau):=x_N(N\tau)$.  For example, $x'_2=x_4^2/(x_4+16)$; also
$x'_3=x_9^3/(x_9^2+9x_9+27)$ and
$x'_4=x_{16}^4/(x_{16}+2)(x_{16}^2+4x_{16}+8)$.  In the same way,
$x'_5=x_{25}^5/(x_{25}^4+5x_{25}^3+15x_{25}^2+25x_{25}+25)$.

If $N=6,7$ then $X_0(N^2)$ is of genus~$1$ but $X_0^+(N^2)=X_0(N^2)/\langle
w_{N^2}\rangle$ is of genus~$0$, and has a Hauptmodul~$t_{N^2}$ that lifts
to a double-valued function on~$X_0(N^2)$.  The following is a
reformulation of Fricke's treatment of this situation, along the lines of
Cohn~\cite{Cohn88}.  A~general point on~$X_0(N^2)$ is determined by the
pair $(t_{N^2},s_{N^2})$, where
\begin{equation}
\label{eq:quartic}
s_{N^2}^2 = {\sf P}_{N^2}(t_{N^2}),
\end{equation}
with ${\sf P}_{N^2}$ a polynomial of degree~$4$ over~$\mathbb{C}$ (with
simple roots).  The projective curve corresponding to this affine quartic
has two infinite points: informally,
$(t_{N^2},s_{N^2})=(\infty,-\infty),(\infty,+\infty)$.  In Fricke's
normalization convention, these are taken to be the distinguished cusps
$\tau={\rm i}\infty,0$, and $t_{N^2}$ is normalized so that the invariants
$j,j_{N^2}$ are asymptotic to $t_{N^2},t_{N^2}^{N^2}$
and~$t_{N^2}^{N^2},t_{N^2}$ as~$t_{N^2}\to\infty$ along the branches
leading to $\tau={\rm i}\infty$ and~$\tau=0$
respectively~\cite[(1.11)]{Cohn88}.  Also, ${\sf P}_{N^2}$~is taken to be
monic.  It then mysteriously turns~out that ${\sf P}_{N^2}$ is
in~$\mathbb{Z}[t_{N^2}]$; this may be fortuitous~\cite[\S\,1]{Cohn91}.

It~will be the case that $x_N(\tau)\,x_N(-1/N\tau)=\kappa_N$ for some
$\kappa_N\in\mathbb{N}$.  So, applying the involution
$w_{N^2}:\tau\mapsto-1/N^2\tau$ to the function $x_N(\tau)$ yields
$\kappa_N/x_N(N\tau)$.  It~follows that a rational function
$F_N(t_{N^2},s_{N^2})$ exists on~$X_0(N^2)$ such that
\begin{subequations}
\begin{align}
x_N(\tau) &= F_N(t_{N^2},s_{N^2}),\\ \kappa_N/x_N(N\tau) &=
F_N(t_{N^2},-s_{N^2}).
\end{align}
\end{subequations}
Any such function is linear in~$s_{N^2}$.  In consequence, the sum and
product of $x_N(\tau)$, $\kappa_N/x_N(N\tau)$ are each rational functions
of~$t_{N^2}$ alone.  That~is, the relation between the sum and product is
of genus~$0$: it~may be uniformized by rational functions.
\makeatletter
\tagsleft@false
\makeatother
\begin{itemize}
\item $N=6$.  $x_6:=72\,\,[2][6]^5/\,[1]^5[3]$, and $x_6$~satisfies
$x_6(\tau)\,x_6(-1/6\tau)=72$.  The $12$-sheeted cover of
$X(1)\cong\mathbb{P}^1(\mathbb{C})_j$ by
$X_0(6)\cong\mathbb{P}^1(\mathbb{C})_{x_6}$ is given by
\begin{align}
\label{eq:12sheeted}
j=\frac{(x_6+6)^3(x_6^3+18\,x_6^2+84\,x_6+24)^3}{x_6\,(x_6+9)^2(x_6+8)^3}.
\end{align}
As elements of the function field of~$X_0(6)$, the functions $x_2,x_3$ may
be expressed rationally as
\begin{align}
\label{eq:x2x3}
x_2=\frac{x_6(x_6+8)^3}{x_6+9},\qquad x_3=\frac{x_6(x_6+9)^2}{x_6+8}.
\end{align}
Combined with the eta-product formulas for $x_2,x_3,x_6$, these imply
\begin{subequations}
\begin{align}
x_6+8&= 8\,\,[2]^9[3]^3/\,[1]^9[6]^3,
\\
x_6+9&= 9\,\,[2]^4[3]^8/\,[1]^8[6]^4.
\end{align}
\end{subequations}
The double cover of $X_0^+(36)\ni t_{36}$ by~$X_0(36)\ni(t_{36},s_{36})$
is given by
\begin{align}
\label{eq:X036eqn}
s_{36}^2 = t_{36}^4 - 8\,t_{36}^3 + 12\,t_{36}^2 - 8\,t_{36} + 4,
\end{align}
from which it follows that as an elliptic curve, $X_0(36)$ is
equianharmonic: it has $j$-invariant zero.  It~is a $6$-sheeted cover
of~$X_0(6)$.  As~elements of the function field of~$X_0(36)$, the functions
$x_6(\tau)$,~$72/x_6(6\tau)$ may be written as
\begin{align}
&\tfrac12(t_{36}-2)\\
&\times\left[(t_{36}^5-10t_{36}^4+28t_{36}^3-26t_{36}^2+20t_{36}+4)\pm(t_{36}^3-6t_{36}^2+6t_{36}-2)s_{36}\right].\nonumber
\end{align}
The formulas
\begin{subequations}
\begin{align}
&x_6(\tau)+72/x_6(6\tau) = (t_{36}-2)\\
&\qquad\qquad\qquad\qquad\qquad\times (t_{36}^5-10\,t_{36}^4+28\,t_{36}^3-26\,t_{36}^2+20\,t_{36}+4),\nonumber\\[\jot]
&x_6(\tau)\cdot[72/x_6(6\tau)] = 72\,t_{36}(t_{36}-2)^2(t_{36}^2-\,t_{36}+1)
\label{eq:x6prod}
\end{align}
\end{subequations}
uniformize the sum--product relation.
\item $N=7$. $x_7:=49\,\,[7]^4/\,[1]^4$, and $x_7$~satisfies
$x_7(\tau)\,x_7(-1/7\tau)=49$.  The $8$-sheeted cover of
$X(1)\cong\mathbb{P}^1(\mathbb{C})_j$ by
$X_0(7)\cong\mathbb{P}^1(\mathbb{C})_{x_7}$ is given by
\begin{align}
j=\frac{(x_7^2+5\,x_7+1)^3(x_7^2+13\,x_7+49)}{x_7}.
\end{align}
The double cover of $X_0^+(49)\ni t_{49}$ by~$X_0(49)\ni(t_{49},s_{49})$ is
given by
\begin{align}
s_{49}^2 = t_{49}^4 - 14\,t_{49}^3 + 63\,t_{49}^2 - 98\,t_{49} + 21,
\end{align}
from which it follows that as an elliptic curve, $X_0(49)$ has
$j$-invariant~$-15^3$.  It~is a $7$-sheeted cover of~$X_0(7)$.  As~elements
of the function field of~$X_0(49)$, the functions
$x_7(\tau)$,~$49/x_7(7\tau)$ may be written as
\begin{align}
&\tfrac12(t_{49}^3-7\,t_{49}^2+14\,t_{49}-7)\\
&\qquad\times\left[(t_{49}^4-14\,t_{49}^3+63\,t_{49}^2-98\,t_{49}+35)\pm(t_{49}^2-7\,t_{49}+7)\,s_{49}\right].\nonumber
\end{align}
The formulas
\begin{subequations}
\begin{align}
&x_7(\tau)+49/x_7(7\tau) = (t_{49}^3-7\,t_{49}^2+14\,t_{49}-7) \\
&\quad\qquad\qquad\qquad\qquad\times(t_{49}^4-14\,t_{49}^3+63\,t_{49}^2-98\,t_{49}+35),\nonumber\\[\jot]
&x_7(\tau)\cdot[49/x_7(7\tau)] = 49\,(t_{49}^3-7\,t_{49}^2+14\,t_{49}-7)^2
\end{align}
\end{subequations}
uniformize the sum--product relation.
\end{itemize}
\makeatletter
\tagsleft@true
\makeatother

\noindent
The above formulas for the $N=6,7$ cases follow from those given by Fricke,
with a good bit of massaging.  For instance, he constructs $X_0(36)$ as a
double cover of the genus-$0$ curve $X_0(18)$, and works~out not a formula
for $x_6$ in~terms of~$(t_{36},s_{36})$, but rather a formula for the
Hauptmodul $x_{18}:=6\,\,[2][3][18]^2/\,[1]^2[6][9]$ of~$X_0(18)$.
However, $x_{6}=x_{18}(x_{18}^2 + 6x_{18} + 12)$, as he computes elsewhere.

It~can be shown that $t_{36}-1=[4][9]\,/\,[1][36]$, relating his
quartic equation~(\ref{eq:X036eqn}) to the equivalent equation
for~$X_0(36)$ derived by Gonz{\'a}lez Rovira~\cite[\S\,4.3]{Gonzalez91}.
Two other notable identities involving~$t_{36}$ were obtained by
Kluit~\cite{Kluit77}, namely
\begin{subequations}
\label{eq:Kluitidents}
\begin{align}
&t_{36}=[2]^3[3][12][18]^3/\,[1]^2[4][6]^2[9][36]^2,\\
&t_{36}-2=[4][6]^8[9]\,/\,[2]^2[3]^3[12]^3[18]^2.
\end{align}
(These eta-products also appear in the tabulation of Ford
et~al.~\cite[Table~4]{Ford94}; see line~36D.)  The further eta-function
identity
\begin{equation}
t_{36}^2-t_{36}+1=[2]^2[3]^4[12]^4[18]^2\,/\,[1]^3[4][6]^4[9][36]^3
\end{equation}
\end{subequations}
follows by substituting (\ref{eq:Kluitidents}ab) into the right side
of~(\ref{eq:x6prod}).  Level~$36$ is in~fact rich in relations among eta
products.  For identities relating eta products that are weight-$1$ modular
forms on the quotient curve~$X_0^+(36)$, see K\"ohler~\cite{Kohler2001}.

\section*{Acknowledgements}

The author would like to thank the referee for valuable comments and
suggestions, which strengthened this article.


\providecommand{\bysame}{\leavevmode\hbox to3em{\hrulefill}\thinspace}
\providecommand{\MR}{\relax\ifhmode\unskip\space\fi MR }
\providecommand{\MRhref}[2]{%
  \href{http://www.ams.org/mathscinet-getitem?mr=#1}{#2}
}
\providecommand{\href}[2]{#2}

\end{document}